\pdfoutput=1
\RequirePackage{ifpdf}
\ifpdf 
\documentclass[pdftex]{sigma}
\else
\documentclass{sigma}
\fi

\usepackage{mathrsfs}

\numberwithin{equation}{section}

\newtheorem{Theorem}{Theorem}[section]
\newtheorem{Lemma}[Theorem]{Lemma}
\newtheorem{Proposition}[Theorem]{Proposition}

\DeclareMathOperator{\trunc}{trunc}
\DeclareMathOperator*{\rep}{Rep}
\DeclareMathOperator*{\diag}{diag}

\begin{document}

\allowdisplaybreaks

\renewcommand{\thefootnote}{$\star$}

\renewcommand{\PaperNumber}{054}

\FirstPageHeading

\ShortArticleName{Extended $T$-System of Type $G_2$}

\ArticleName{Extended $\boldsymbol{T}$-System of Type $\boldsymbol{G_2}$\footnote{This paper is
a~contribution to the Special Issue in honor of Anatol Kirillov and Tetsuji Miwa. The full collection is available at
\href{http://www.emis.de/journals/SIGMA/InfiniteAnalysis2013.html}
{http://www.emis.de/journals/SIGMA/InfiniteAnalysis2013.html}}}

\Author{Jian-Rong LI~$^\dag$ and Evgeny MUKHIN~$^\ddag$}

\AuthorNameForHeading{J.R.~Li and E.~Mukhin}

\Address{$^\dag$~Department of Mathematics, Lanzhou University, Lanzhou 730000, P.R.~China}
\EmailD{\href{mailto:lijr@lzu.edu.cn}{lijr@lzu.edu.cn}, \href{mailto:lijr07@gmail.com}{lijr07@gmail.com}}
\URLaddressD{\url{http://scholar.google.com/citations?user=v_0AZ7oAAAAJ&hl=en}}

\Address{$^\ddag$~Department of Mathematical Sciences, Indiana University - Purdue University Indianapolis,\\
\hphantom{$^\ddag$}~402 North Blackford St, Indianapolis, IN 46202-3216, USA}
\EmailD{\href{mailto:mukhin@math.iupui.edu}{mukhin@math.iupui.edu}}
\URLaddressD{\url{http://www.math.iupui.edu/~mukhin/}}

\ArticleDates{Received April 03, 2013, in f\/inal form August 16, 2013; Published online August 22, 2013}

\Abstract{We prove a~family of $3$-term relations in the Grothendieck ring of the category of
f\/inite-dimensional modules over the af\/f\/ine quantum algebra of type $G_2$ extending the celebrated
$T$-system relations of type $G_2$.
We show that these relations can be used to compute classes of certain irreducible modules, including
classes of all minimal af\/f\/inizations of type $G_2$.
We use this result to obtain explicit formulas for dimensions of all participating modules.}


\Classification{17B37; 81R50; 82B23}

\renewcommand{\thefootnote}{\arabic{footnote}}
\setcounter{footnote}{0}

\section{Introduction}

Kirillov--Reshetikhin modules are simplest examples of irreducible f\/inite-dimensional modules over
quantum af\/f\/ine algebras, and the $T$-system is a~famous family of short exact sequences of tensor
products of Kirillov--Reshetikhin modules, see~\cite{Her06, KR90,KNS94,Nak03}.
There are numerous applications of the $T$-systems in representation theory, combinatorics and integrable
systems, see the survey~\cite{KNS11}.

Minimal af\/f\/inizations of quantum af\/f\/ine algebras form an important family of irreducible modules
which contains the Kirillov--Reshetikhin modules, see~\cite{CP95b}.
A procedure to extend the $T$-system to a~larger set of relations to include the minimal af\/f\/inization was
described in~\cite{MY11b}, where it was conjectured to work in all types.
In~\cite{MY11b} this procedure was carried out in types $A$ and $B$.
In this paper, we show the existence of the extended $T$-system for type $G_2$.

We work with the quantum af\/f\/ine algebra $U_q\hat{\mathfrak{g}}$ of type $G_2$.
The irreducible f\/inite-dimensional modules of quantum af\/f\/ine algebras are parameterized by the
highest $l$-weights or Drinfeld polynomials.
Let $\mathcal{T}$ be an irreducible $U_q\hat{\mathfrak{g}}$-module such that zeros of all Drinfeld
polynomials belong to a~lattice $aq^\mathbb{Z}$ for some $a\in\mathbb{C}^{\times}$.
Following~\cite{MY11b}, we def\/ine the left, right, and bottom modules, denoted by $\mathcal{L}$,
$\mathcal{R}$, $\mathcal{B}$ respectively.
The Drinfeld polynomials of left, right, and bottom modules are obtained by stripping the rightmost,
leftmost, and both left- and rightmost zeros of the union of zeros of the Drinfeld polynomials of the top
module $\mathcal{T}$.

Then the relations of the extended $T$-system have the form $[\mathcal{L}] [\mathcal{R}]=[\mathcal{T}]
[\mathcal{B}] + [\mathcal{S}]$, where $[ \cdot ]$ denotes the equivalence class of
a~$U_q\hat{\mathfrak{g}}$-module in the Grothendieck ring of the category of f\/inite-dimensional
representations of $U_q\hat{\mathfrak{g}}$.
Moreover, in all cases the modules $\mathcal{T}\otimes \mathcal{B}$ and $\mathcal{S}$ are irreducible.

We start with minimal af\/f\/inizations as the top modules $\mathcal{T}$, then the left, right and bottom
modules are minimal af\/f\/inizations as well.
We compute $S$ and decompose it as a~product of irreducible modules which we call sources.
It turns out that the sources are not always minimal af\/f\/inizations.
Therefore, we follow up with taking the sources as top modules and compute new left, right, bottom modules,
and sources.
Then we use all new modules obtained on a~previous step as top modules and so on.

We end up with several families of modules which we denote by $\mathcal{B}_{k, \ell}^{(s)}$,
$\mathcal{C}_{k, \ell}^{(s)}$, $\mathcal{D}_{k, \ell}^{(s)}$, $\mathcal{E}_{k, \ell}^{(s)}$,
$\mathcal{F}_{k, \ell}^{(s)}$, $\tilde{\mathcal{B}}_{k, \ell}^{(s)}$, $\tilde{\mathcal{C}}_{k,
\ell}^{(s)}$, $\tilde{\mathcal{D}}_{k, \ell}^{(s)}$, $\tilde{\mathcal{E}}_{k, \ell}^{(s)}$,
$\tilde{\mathcal{F}}_{k, \ell}^{(s)}$, where $s\in \mathbb{Z}$, $k,\ell\in \mathbb{Z}_{\geq 0}$.
This is the minimal set of modules which contains all minimal af\/f\/inizations \big(these are modules
$\mathcal{B}_{k, \ell}^{(s)}$, $\tilde{\mathcal{B}}_{k, \ell}^{(s)}$\big) and which is closed under our set of
relations.
Namely, if any of the above modules is chosen as a~top module then the left, right, bottom modules and all
sources belong to this set as well, see Theorems~\ref{extended T-system},~\ref{extended T-system for tilde}.
The spirit of the proof Theorems~\ref{extended T-system},~\ref{extended T-system for tilde} follows the
works in~\cite{Her06,MY11b, Nak04}.

We show that the extended $T$-system allows us to compute the modules $\mathcal{B}_{k, \ell}^{(s)}$,
$\mathcal{C}_{k, \ell}^{(s)}$, $\mathcal{D}_{k, \ell}^{(s)}$, $\mathcal{E}_{k, \ell}^{(s)}$,
$\mathcal{F}_{k, \ell}^{(s)}$ recursively in terms of fundamental modules, see Proposition~\ref{compute}.
We use this to compute the dimensions of all participating modules, in particular, we give explicit
formulas for dimensions of all minimal af\/f\/inizations of type $G_2$, see Theorem~\ref{dimensions}.
We hope further, that one can use use the extended $T$-system to obtain the decomposition of all
participating modules as the $U_q{\mathfrak{g}}$-modules.

Let us point out some similarities and dif\/ferences with types $A$ and $B$.
The type $A$, the extended $T$-system is closed within the class of minimal af\/f\/inizations, meaning that
all sources are minimal af\/f\/inizations as well.
In type $B$, the extended $T$-system is not closed within the class of minimal af\/f\/inizations, but it is
closed in the class of so called snake modules, see~\cite{MY11b}.
For the proofs and computations it is important that all modules participating in extended $T$-systems of
types $A$ and $B$ are thin and special, moreover their $q$-characters are known explicitly in terms of skew
Young tableaux in type $A$, and in terms of path models in type $B$, see~\cite{Che87, MY11b, MY11a,NT98}.

In general the modules of the extended $T$-system of type $G_2$ are not thin and at the moment there is no
combinatorial description of their $q$-characters.
However, all modules turn out to be either special or anti-special.
Therefore we are able to use the FM algorithm, see~\cite{FM01}, to compute the suf\/f\/icient information
about $q$-characters in order to complete the proofs.
Note, that a~priori it is was not obvious that the extended $T$-system will be closed within special or
anti-special modules.
Moreover, since the $q$-characters of $G_2$ modules are not known explicitly, the property of being special
or anti-special had to be established in each case, see Theorems~\ref{special},~\ref{anti-special}.

Note that in general the minimal af\/f\/inizations of types $C$, $D$, $E$, $F$ are neither special nor
anti-special, therefore the methods of this paper cannot be applied in those cases.

There is a~remarkable conjecture on the cluster algebra relations in the category of f\/inite-dimensional
representations of quantum af\/f\/ine algebras of type $A$, $D$, $E$, see~\cite{HL10}.
Taking into account the work of~\cite{IIKKN10a,IIKKN10b}, one could expect that the conjecture
of~\cite{HL10} can be formulated for other types as well, in particular for type~$G_2$.
We expect that the extended $T$-system is a~part of cluster algebra relations.

The paper is organized as follows.
In Section~\ref{background}, we give some background material.
In Section~\ref{main results}, we def\/ine the modules $\mathcal{B}_{k, \ell}^{(s)}$,
$\mathcal{C}_{k, \ell}^{(s)}$, $\mathcal{D}_{k, \ell}^{(s)}$, $\mathcal{E}_{k, \ell}^{(s)}$,
$\mathcal{F}_{k, \ell}^{(s)}$ and state our main result, Theorem~\ref{extended T-system}.
In Section~\ref{proof special}, we prove that the modules $\mathcal{B}_{k, \ell}^{(s)}$,
$\mathcal{C}_{k, \ell}^{(s)}$, $\mathcal{D}_{k, \ell}^{(s)}$, $\mathcal{E}_{k, \ell}^{(s)}$,
$\mathcal{F}_{k, \ell}^{(s)}$ are special.
In Section~\ref{proof system}, we prove Theorem~\ref{extended T-system}.
In Section~\ref{proof irreducible}, we prove that the module $\mathcal{T} \otimes \mathcal{B}$ is
irreducible for each relation in the extended $T$-system.
In Section~\ref{The second part of the extended T-system}, we deduce the extended $T$-system for the
modules $\tilde{\mathcal{B}}_{k, \ell}^{(s)}$, $\tilde{\mathcal{C}}_{k, \ell}^{(s)}$,
$\tilde{\mathcal{D}}_{k, \ell}^{(s)}$, $\tilde{\mathcal{E}}_{k, \ell}^{(s)}$, $\tilde{\mathcal{F}}_{k,
\ell}^{(s)}$.
In Section~\ref{dimensions section}, we compute the dimensions of the modules in the extended
$T$-systems.

\section{Background}
\label{background}

\subsection{Cartan data}

Let $\mathfrak{g}$ be a~complex simple Lie algebra of type $G_2$ and
$\mathfrak{h}$ a~Cartan subalgebra of $\mathfrak{g}$.
Let $I=\{1, 2\}$.
We choose simple roots $\alpha_1$, $\alpha_2$ and scalar product $(\cdot, \cdot)$ such that
\begin{gather*}
(\alpha_1,\alpha_1)=2,
\qquad
(\alpha_1,\alpha_2)=-3,
\qquad
(\alpha_2,\alpha_2)=6.
\end{gather*}
Let $\{\alpha_1^{\vee}, \alpha_2^{\vee}\}$ and $\{\omega_1, \omega_2\}$ be the sets of simple coroots and
fundamental weights respectively.
Let $C=(C_{ij})_{i, j\in I}$ denote the Cartan matrix, where $C_{ij}=\frac{2 ( \alpha_i, \alpha_j ) }{(
\alpha_i, \alpha_i )}$.
Let $r_1=1$, $r_2=3$, $D=\diag (r_1, r_2)$ and $B=DC$.
Then
\begin{gather*}
C=\left(
\begin{matrix}2&-3
\\
-1&2
\end{matrix}
\right),
\qquad
B=\left(
\begin{matrix}2&-3
\\
-3&6
\end{matrix}
\right).
\end{gather*}

Let $Q$ (resp.\
{}$Q^+$) and $P$ (resp.\
{}$P^+$) denote the $\mathbb{Z}$-span (resp.\
{}$\mathbb{Z}_{\geq 0}$-span) of the simple roots and fundamental weights respectively.
Let $\leq$ be the partial order on $P$ in which $\lambda \leq \lambda'$ if and only if $\lambda' - \lambda
\in Q^+$.

Let $\hat{\mathfrak{g}}$ denote the untwisted af\/f\/ine algebra corresponding to~$\mathfrak{g}$.
Fix a~$q\in \mathbb{C}^{\times}$, not a~root of unity.
Let $q_i=q^{r_i}$, $i=1, 2$.
Def\/ine the $q$-numbers, $q$-factorial and $q$-binomial:
\begin{gather*}
[n]_q:=\frac{q^n-q^{-n}}{q-q^{-1}},
\qquad
[n]_q!:=[n]_q[n-1]_q\cdots[1]_q,
\qquad
{n\brack m}_q:=\frac{[n]_q!}
{[n-m]_q![m]_q!}.
\end{gather*}

\subsection{Quantum af\/f\/ine algebra}

The quantum af\/f\/ine algebra $U_q\hat{\mathfrak{g}}$ in Drinfeld's new realization, see~\cite{Dri88}, is
generated by $x_{i, n}^{\pm}$ ($i\in I$, $n\in \mathbb{Z}$), $k_i^{\pm 1}$ $(i\in I)$, $h_{i, n}$ ($i\in I$,
$n\in \mathbb{Z}\backslash \{0\}$) and central elements $c^{\pm 1/2}$, subject to the following relations:
\begin{gather*}
k_i k_j=k_j k_i,
\qquad
k_i h_{j,n}=h_{j,n}k_i,
\qquad
k_ik_i^{-1}=k_i^{-1}k_i=1,
\qquad
k_i x_{j,n}^{\pm}k_i^{-1}=q^{\pm B_{ij}}x_{j,n}^{\pm},
\\
\big[h_{i, n}, x_{j, m}^{\pm}\big]
=\pm\frac{1}{n}[nB_{ij}]_q c^{\mp|n|/2}x_{j,n+m}^{\pm},
\\
x_{i,n+1}^{\pm}x_{j,m}^{\pm}-q^{\pm B_{ij}}x_{j,m}^{\pm}x_{i,n+1}^{\pm}=q^{\pm B_{ij}}x_{i,n}^{\pm}
x_{j,m+1}^{\pm}-x_{j,m+1}^{\pm}x_{i,n}^{\pm},
\\
[h_{i, n}, h_{j, m}]
=\delta_{n,-m}\frac{1}{n}[nB_{ij}]_q\frac{c^n-c^{-n}}{q-q^{-1}},
\\
[x_{i, n}^{+}, x_{j, m}^{-}]
=\delta_{ij}\frac{c^{(n-m)/2}\phi_{i,n+m}^{+}-c^{-(n-m)/2}\phi_{i,n+m}^{-}}{q_i-q_i^{-1}},
\\
\sum_{\pi\in\Sigma_s}\sum_{k=0}^{s}(-1)^k{s\brack k}_{q_i}x_{i,n_{\pi(1)}}^{\pm}\cdots x_{i,n_{\pi(k)}}
^{\pm}x_{j,m}^{\pm}x_{i,n_{\pi(k+1)}}^{\pm}\cdots x_{i,n_{\pi(s)}}^{\pm}=0,
\qquad
s=1-C_{ij},
\end{gather*}
for all sequences of integers $n_1,\ldots,n_s$, and $i \neq j$, where $\Sigma_{s}$ is the symmetric
groups on $s$ letters, $\phi_{i, n}^{\pm}=0$ ($n < 0$) and $\phi_{i, n}^{\pm}$'s ($n \geq 0$) are
determined by the formula
\begin{gather}
\phi_i^{\pm}(u):=\sum_{n=0}^{\infty}\phi_{i,\pm n}^{\pm}u^{\pm n}=k_i^{\pm1}\exp\left(\pm\big(q-q^{-1}
\big)\sum_{m=1}^{\infty}h_{i,\pm m}u^{\pm m}\right).
\end{gather}
There exist a~coproduct, counit and antipode making $U_q\hat{\mathfrak{g}}$ into a~Hopf algebra.

The quantum af\/f\/ine algebra $U_q\hat{\mathfrak{g}}$ contains two standard quantum af\/f\/ine algebras of
type~$A_1$.
The f\/irst one is $U_{q_1}(\hat{\mathfrak{sl}_2})$ generated by $x_{1, n}^{\pm}$ ($n\in \mathbb{Z}$),
$k_1^{\pm 1}$, $h_{1, n}$ ($n\in \mathbb{Z}\backslash \{0\}$) and central elements~$c^{\pm 1/2}$.
The second one is $U_{q_2}(\hat{\mathfrak{sl}_2})$ generated by $x_{2, n}^{\pm}$ ($n\in \mathbb{Z}$),
$k_2^{\pm 1}$, $h_{2, n}$ ($n\in \mathbb{Z}\backslash \{0\}$) and central elements~$c^{\pm 1/2}$.

The subalgebra of $U_q\hat{\mathfrak{g}}$ generated by $(k_i^{\pm})_{i\in I}$, $(x_{i, 0}^{\pm})_{i\in I}$ is
a~Hopf subalgebra of $U_q\hat{\mathfrak{g}}$ and is isomorphic as a~Hopf algebra to $U_q\mathfrak{g}$, the
quantized enveloping algebra of $\mathfrak{g}$.
In this way, $U_q\hat{\mathfrak{g}}$-modules restrict to $U_q\mathfrak{g}$-modules.

\subsection[Finite-dimensional representations and $q$-characters]
{Finite-dimensional representations and $\boldsymbol{q}$-characters}

In this section, we recall the standard
facts about f\/inite-dimensional representations of $U_q\hat{\mathfrak{g}}$ and $q$-characters of these
representations, see~\cite{CP94,CP95a,FM01,FR98,MY11b}.

A representation $V$ of $U_q\hat{\mathfrak{g}}$ is of type~$1$ if $c^{\pm 1/2}$ acts as the identity on $V$
and
\begin{gather}
\label{decomposition}
V=\bigoplus_{\lambda\in P}V_{\lambda},
\qquad
V_{\lambda}=\big\{v\in V:k_i v=q^{(\alpha_i,\lambda)}v\big\}.
\end{gather}
In the following, all representations will be assumed to be f\/inite-dimensional and of type~$1$ without
further comment.
The decomposition~\eqref{decomposition} of a~f\/inite-dimensional representation $V$ into its
$U_q\mathfrak{g}$-weight spaces can be ref\/ined by decomposing it into the Jordan subspaces of the
mutually commuting operators $\phi_{i, \pm r}^{\pm}$, see~\cite{FR98}:
\begin{gather}
V=\bigoplus_{\gamma}V_{\gamma},
\qquad
\gamma=\big(\gamma_{i,\pm r}^{\pm}\big)_{i\in I,\;r\in\mathbb{Z}_{\geq0}}
,
\qquad
\gamma_{i,\pm r}^{\pm}\in\mathbb{C},
\end{gather}
where
\begin{gather*}
V_{\gamma}=\left\{v\in V:\exists\, k\in\mathbb{N},\forall\, i\in I,m\geq0,\big(\phi_{i,\pm m}^{\pm}-\gamma_{i,\pm m}
^{\pm}\big)^{k}v=0\right\}.
\end{gather*}
If $\dim(V_{\gamma})>0$, then $\gamma$ is called an \textit{$l$-weight} of $V$.
For every f\/inite dimensional representation of~$U_q\hat{\mathfrak{g}}$, the $l$-weights are known,
see~\cite{FR98}, to be of the form
\begin{gather}
\gamma_i^{\pm}(u):=\sum_{r=0}^{\infty}\gamma_{i,\pm r}^{\pm}u^{\pm r}=q_i^{\deg Q_i-\deg R_i}
\frac{Q_i(uq_i^{-1})R_i(uq_i)}{Q_{i}(uq_i)R_{i}(uq_i^{-1})},
\label{gamma}
\end{gather}
where the right hand side is to be treated as a~formal series in positive (resp.\
negative) integer powers of $u$, and $Q_i$, $R_i$ are polynomials of the form
\begin{gather}
Q_i(u)=\prod_{a\in\mathbb{C}^{\times}}(1-ua)^{w_{i,a}},
\qquad
R_{i}(u)=\prod_{a\in\mathbb{C}^{\times}}
(1-ua)^{x_{i,a}},
\label{QR}
\end{gather}
for some $w_{i, a}$, $x_{i, a} \in \mathbb{Z}_{\geq 0}$, $i\in I$, $a\in \mathbb{C}^{\times}$.
Let $\mathcal{P}$ denote the free abelian multiplicative group of monomials in inf\/initely many formal
variables $(Y_{i, a})_{i\in I,\;a\in\mathbb{C}^{\times}}$.
There is a~bijection~$\gamma$ from $\mathcal{P}$ to the set of $l$-weights of f\/inite-dimensional modules
such that for the monomial $m=\prod_{i\in I,\;a\in\mathbb{C}^{\times}} Y_{i, a}^{w_{i, a}-x_{i,
a}}$, the $l$-weight $\gamma(m)$ is given by~\eqref{gamma},~\eqref{QR}.

Let $\mathbb{Z}\mathcal{P} = \mathbb{Z}\big[Y_{i, a}^{\pm 1}\big]_{i\in I,\;a\in\mathbb{C}^{\times}}$ be the group
ring of $\mathcal{P}$.
For $\chi \in \mathbb{Z}\mathcal{P}$, we write $m\in \mathcal{P}$ if the coef\/f\/icient of $m$ in $\chi$
is non-zero.

The $q$-character of a~$U_q\hat{\mathfrak{g}}$-module $V$ is given by
\begin{gather*}
\chi_q(V)=\sum_{m\in\mathcal{P}}\dim(V_{m})m\in\mathbb{Z}\mathcal{P},
\end{gather*}
where $V_m = V_{\gamma(m)}$.

Let $\rep(U_q\hat{\mathfrak{g}})$ be the Grothendieck ring of f\/inite-dimensional representations of
$U_q\hat{\mathfrak{g}}$ and $[V]\in \rep(U_q\hat{\mathfrak{g}})$ the class of a~f\/inite-dimensional
$U_q\hat{\mathfrak{g}}$-module $V$.
The $q$-character map def\/ines an injective ring homomorphism, see~\cite{FR98},
\begin{gather*}
\chi_q: \ \rep(U_q\hat{\mathfrak{g}})\to\mathbb{Z}\mathcal{P}.
\end{gather*}

For any f\/inite-dimensional representation $V$ of $U_q\hat{\mathfrak{g}}$, denote by $\mathscr{M}(V)$ the
set of all monomials in $\chi_q(V)$.
For each $j\in I$, a~monomial $m=\prod_{i\in I,\;a\in\mathbb{C}^{\times}} Y_{i, a}^{u_{i, a}}$,
where $u_{i, a}$ are some integers, is said to be \textit{$j$-dominant} (resp.\
\textit{$j$-anti-dominant}) if and only if $u_{j, a} \geq 0$ (resp.\
{}$u_{j, a} \leq 0$) for all $a\in \mathbb{C}^{\times}$.
A monomial is called \textit{dominant} (resp.\
\textit{anti-dominant}) if and only if it is $j$-dominant (resp.\
{}$j$-anti-dominant) for all $j\in I$.
Let $\mathcal{P}^+ \subset \mathcal{P}$ denote the set of all dominant monomials.

Let $V$ be a~representation of $U_q\hat{\mathfrak{g}}$ and $m\in \mathscr{M}(V)$ a~monomial.
A non-zero vector $v\in V_m$ is called a~\textit{highest $l$-weight vector} with \textit{highest
$l$-weight} $\gamma(m)$ if
\begin{gather*}
x_{i,r}^{+}\cdot v=0,
\qquad
\phi_{i,\pm t}^{\pm}\cdot v=\gamma(m)_{i,\pm t}^{\pm}
v,
\qquad
\forall\, i\in I,\;r\in\mathbb{Z},\;t\in\mathbb{Z}_{\geq0}.
\end{gather*}
The module $V$ is called a~\textit{highest $l$-weight representation} if $V=U_q\hat{\mathfrak{g}}\cdot v$
for some highest $l$-weight vector $v\in V$.

It is known, see~\cite{CP94,CP95a}, that for each $m_+\in \mathcal{P}^{+}$ there is a~unique
f\/inite-dimensional irreducible representation, denoted $L(m_+)$, of $U_q\hat{\mathfrak{g}}$ that is
highest $l$-weight representation with highest $l$-weight $\gamma(m_+)$, and moreover every
f\/inite-dimensional irreducible $U_q\hat{\mathfrak{g}}$-module is of this form for some $m_+ \in
\mathcal{P}^+$.
Also, if $m_+$, $m_+' \in \mathcal{P}^{+}$ and $m_{+}\neq m_{+}'$, then $L(m_+) \not\cong L(m_+')$.
For $m_+ \in \mathcal{P}^+$, we use $\chi_q(m_+)$ to denote $\chi_q(L(m_+))$.

The following lemma is well-known.
\begin{Lemma}
Let $m_1$, $m_2$ be two monomials.
Then $L(m_1m_2)$ is a~sub-quotient of $L(m_1) \otimes L(m_2)$.
In particular, $\mathscr{M}(L(m_1m_2)) \subseteq \mathscr{M}(L(m_1))\mathscr{M}(L(m_2))$.
\end{Lemma}

For $b\in \mathbb{C}^{\times}$, def\/ine the shift of spectral parameter map $\tau_b: \mathbb{Z}\mathcal{P}
\to \mathbb{Z}\mathcal{P}$ to be a~homomorphism of rings sending $Y_{i, a}^{\pm 1}$ to $Y_{i, ab}^{\pm 1}$.
Let $m_1$, $m_2 \in \mathcal{P}^+$.
If $\tau_b(m_1) = m_2$, then
\begin{gather}
\label{shift}
\tau_b\chi_q(m_1)=\chi_q(m_2).
\end{gather}

Let $m_{+}$ be a~dominant $l$-weight.
We call the polynomial $\chi_q(m_{+})$ special if it contains exactly one dominant monomial.

A f\/inite-dimensional $U_q\hat{\mathfrak{g}}$-module $V$ is said to be \textit{special} if and only if
$\mathscr{M}(V)$ contains exactly one dominant monomial.
It is called \textit{anti-special} if and only if $\mathscr{M}(V)$ contains exactly one anti-dominant
monomial.
It is called \textit{thin} if and only if no $l$-weight space of $V$ has dimension greater than~$1$.
We also call a~polynomial in $\Bbb{Z}\mathcal P$ {\it special}, {\it antispecial}, or {\it thin} if this
polynomial contains a~unique dominant monomial, a~unique anti-dominant monomial, or if all coef\/f\/icients
are zero and one respectively.
A f\/inite-dimensional $U_q\hat{\mathfrak{g}}$-module is said to be \textit{prime} if and only if it is not
isomorphic to a~tensor product of two non-trivial $U_q\hat{\mathfrak{g}}$-modules, see~\cite{CP97}.
Clearly, if a~module is special or anti-special, then it is irreducible.

Def\/ine $A_{i, a} \in \mathcal{P}$, $i\in I$, $a\in \mathbb{C}^{\times}$, by
\begin{gather*}
A_{1,a}=Y_{1,aq}Y_{1,aq^{-1}}Y_{2,a}^{-1},
\qquad
A_{2,a}=Y_{2,aq^{3}}Y_{2,aq^{-3}}Y_{1,aq^{-2}}^{-1}Y_{1,a}^{-1}Y_{1,aq^{2}}^{-1}.
\end{gather*}
Let $\mathcal{Q}$ be the subgroup of $\mathcal{P}$ generated by $A_{i, a}$, $i\in I$, $a\in\mathbb{C}^{\times}$.
Let $\mathcal{Q}^{\pm}$ be the monoids generated by $A_{i, a}^{\pm 1}$, $i\in I$, $a\in\mathbb{C}^{\times}$.
There is a~partial order $\leq$ on $\mathcal{P}$ in which
\begin{gather}
m\leq m' \ \  \text{if and only if} \ \ m'm^{-1}\in\mathcal{Q}^{+}.
\label{partial order of monomials}
\end{gather}
For all $m_+ \in \mathcal{P}^+$, $\mathscr{M}(L(m_+)) \subset m_+\mathcal{Q}^{-}$, see~\cite{FM01}.

A monomial $m$ is called \textit{right negative} if and only if $\forall\, a~\in \mathbb{C}^{\times}$ and
$\forall\, i\in I$, we have the following property: if the power of $Y_{i, a}$ is non-zero and the power of
$Y_{j, aq^k}$ is zero for all $j \in I$, $k\in \mathbb{Z}_{>0}$, then the power of $Y_{i, a}$ is negative.
For $i\in I$, $a\in \mathbb{C}^{\times}$, $A_{i,a}^{-1}$ is right-negative.
A~product of right-negative monomials is right-negative.
If~$m$ is right-negative, then $m'\leq m$ implies that~$m'$ is right-negative.

\subsection[Minimal af\/f\/inizations of $U_q\mathfrak{g}$-modules]
{Minimal af\/f\/inizations of $\boldsymbol{U_q\mathfrak{g}}$-modules}

Let $\lambda = k\omega_1 + \ell \omega_2$.
A simple $U_q\hat{\mathfrak{g}}$-module $L(m_+)$ is called a~\textit{minimal affinization} of $V(\lambda)$
if and only if $m_+$ is one of the following monomials
\begin{gather*}
\left(\prod_{i=0}^{\ell-1}Y_{2,aq^{6i}}\right)\left(\prod_{i=0}^{k-1}Y_{1,aq^{6\ell+2i+1}}\right),
\qquad
\left(\prod_{i=0}^{k-1}Y_{1,aq^{2i}}\right)\left(\prod_{i=0}^{\ell-1}Y_{2,aq^{2k+6i+5}}\right),
\end{gather*}
for some $a\in \mathbb{C}^{\times}$, see~\cite{CP95b}.
In particular, when $k=0$ or $\ell =0$, the minimal af\/f\/inization $L(m_+)$ is called
a~\textit{Kirillov--Reshetikhin module}.

Let $L(m_+)$ be a~Kirillov--Reshetikhin module.
It is shown in~\cite{Her06} that any non-highest monomial in $\mathscr{M}(L(m_+))$ is right-negative and in
particular $L(m_+)$ is special.

\subsection[$q$-characters of $U_q\hat{\mathfrak{sl}}_2$-modules and the FM algorithm]
{$\boldsymbol{q}$-characters of $\boldsymbol{U_q\hat{\mathfrak{sl}}_2}$-modules and the FM algorithm}

The $q$-characters of $U_q\hat{\mathfrak{sl}}_2$-modules are well-understood, see~\cite{CP91,FR98}.
We recall the results here.

Let $W_{k}^{(a)}$ be the irreducible representation $U_q\hat{\mathfrak{sl}}_2$ with highest weight monomial
\begin{gather*}
X_{k}^{(a)}=\prod_{i=0}^{k-1}Y_{aq^{k-2i-1}},
\end{gather*}
where $Y_a=Y_{1, a}$.
Then the $q$-character of $W_{k}^{(a)}$ is given by
\begin{gather*}
\chi_q\big(W_{k}^{(a)}\big)=X_{k}^{(a)}\sum_{i=0}^{k}\prod_{j=0}^{i-1}A_{aq^{k-2j}}^{-1},
\end{gather*}
where $A_a=Y_{aq^{-1}}Y_{aq}$.

For $a\in \mathbb{C}^{\times}$, $k\in \mathbb{Z}_{\geq 1}$,
the set $\Sigma_{k}^{(a)}=\{aq^{k-2i-1}\}_{i=0,\ldots,k-1}$ is called a~\textit{string}.
Two strings $\Sigma_{k}^{(a)}$ and $\Sigma_{k'}^{(a')}$ are said to be in \textit{general position} if the
union $\Sigma_{k}^{(a)} \cup \Sigma_{k'}^{(a')}$ is not a~string or $\Sigma_{k}^{(a)} \subset
\Sigma_{k'}^{(a')}$ or $\Sigma_{k'}^{(a')} \subset \Sigma_{k}^{(a)}$.

Denote by $L(m_+)$ the irreducible $U_q\hat{\mathfrak{sl}}_2$-module with highest weight monomial $m_+$.
Let \mbox{$m_{+} \neq 1$} and $\in \mathbb{Z}[Y_a]_{a\in \mathbb{C}^{\times}}$ be a~dominant monomial.
Then $m_+$ can be uniquely (up to permutation) written in the form
\begin{gather*}
m_+=\prod_{i=1}^{s}\left(\prod_{b\in\Sigma_{k_i}^{(a_i)}}Y_{b}\right),
\end{gather*}
where $s$ is an integer, $\Sigma_{k_i}^{(a_i)}$, $i=1,\ldots,s$, are strings which are pairwise in general
position and
\begin{gather*}
L(m_+)=\bigotimes_{i=1}^s W_{k_i}^{(a_i)},
\qquad
\chi_q(m_+)=\prod_{i=1}^s\chi_q\Big(W_{k_i}^{(a_i)}\Big).
\end{gather*}

For $j\in I$, let
\begin{gather*}
\beta_j:\ \mathbb{Z}\big[Y_{i,a}^{\pm1}\big]_{i\in I;a\in\mathbb{C}^{\times}}
\to\mathbb{Z}\big[Y_{a}^{\pm1}\big]_{a\in\mathbb{C}^{\times}}
\end{gather*}
be the ring homomorphism which sends, for all $a\in \mathbb{C}^{\times}$, $Y_{k, a} \mapsto 1$ for $k\neq
j$ and $Y_{j, a} \mapsto Y_{a}$.

Let $V$ be a~$U_q\hat{\mathfrak{g}}$-module.
Then $\beta_i(\chi_q(V))$, $i=1, 2$, is the $q$-character of $V$ considered as
a~$U_{q_i}(\hat{\mathfrak{sl}_2})$-module.

In some situation, we can use the $q$-characters of $U_q\hat{\mathfrak{sl}}_2$-modules to compute the
$q$-characters of $U_q\hat{\mathfrak{g}}$-modules for arbitrary $\mathfrak{g}$, see~\cite[Section~5]{FM01}.
The corresponding algorithm is called the FM algorithm.
The FM algorithm recursively computes the minimal possible $q$-character which contains $m_+$ and is
consistent when restricted to $U_{q_i}(\hat{\mathfrak{sl}_2})$, $i=1, 2$.

Although the FM algorithm does not give the $q$-character of a~$U_q\hat{\mathfrak{g}}$-module in general,
the FM algorithm works for a~large family of $U_q\hat{\mathfrak{g}}$-modules.
For example, if a~module $L(m_+)$ is special, then the FM algorithm applied to $m_+$, produces the correct
$q$-character $\chi_q(m_+)$, see~\cite{FM01}.

\subsection[Truncated $q$-characters]{Truncated $\boldsymbol{q}$-characters}

We use the truncated $q$-characters \cite{HL10,MY11b}.
Given a~set of monomials $\mathcal{R} \subset \mathcal{P}$, let $\mathbb{Z}\mathcal{R} \subset
\mathbb{Z}\mathcal{P}$ denote the $\mathbb{Z}$-module of formal linear combinations of elements of
$\mathcal{R}$ with integer coef\/f\/icients.
Def\/ine
\begin{gather*}
\trunc_{\mathcal{R}}: \ \mathcal{P}\to\mathcal{R};
\qquad
m\mapsto
\begin{cases}
m&\text{if}\quad m\in\mathcal{R},
\\
0&\text{if}\quad m\not\in\mathcal{R},
\end{cases}
\end{gather*}
and extend $\trunc_{\mathcal{R}}$ as a~$\mathbb{Z}$-module map $\mathbb{Z}\mathcal{P} \to
\mathbb{Z}\mathcal{R}$.

Given a~subset $U\subset I \times \mathbb{C}^{\times}$, let $\mathcal{Q}_U$ be the subgroups of
$\mathcal{Q}$ generated by $A_{i,a}$ with $(i,a)\in U$.
Let $\mathcal{Q}_U^{\pm}$ be the monoid generated by $A_{i,a}^{\pm 1}$ with $(i,a)\in U$.
We call $\trunc_{m_{+}\mathcal{Q}_U^{-}}   \chi_q(m_{+})$ \textit{the $q$-character of $L(m_{+})$
truncated to $U$}.

If $U=I\times \mathbb{C}^{\times}$, then $\trunc_{m_{+}\mathcal{Q}_U^{-}}   \chi_q(m_{+})$ is the
ordinary $q$-character of $L(m_{+})$.

The main idea of using the truncated $q$-characters is the following.
Given $m^+$, one chooses $\mathcal R$ in such a~way that the dropped monomials are all right-negative and
the truncated $q$-character is much smaller than the full $q$-character.
The advantage is that the truncated $q$-character is much easier to compute and to describe in
a~combinatorial way.
At the same time, if the truncating set $\mathcal R$ can be used for both $m_1^+$ and $m_2^+$, then the
same $\mathcal R$ works for the tensor product $L(m_1^+)\otimes L(m_2^+)$.
Moreover, the product of truncated characters of $L(m_1^+)$ and $L(m_2^+)$ contains all dominant monomials
of the tensor product $L(m_1^+)\otimes L(m_2^+)$ and can be used to f\/ind the decomposition of it into
irreducible components in the Grothendieck ring.
We compute the truncated $q$-characters using the following theorem.
\begin{Theorem}[\protect{\cite[Theorem 2.1]{MY11b}}]
\label{truncated}
Let $U\subset I \times \mathbb{C}^{\times}$ and $m_{+} \in \mathcal{P}^+$.
Suppose that $\mathcal{M} \subset \mathcal{P}$ is a~finite set of distinct monomials such that
\begin{enumerate}\itemsep=0pt
\item[$(i)$] $\mathcal{M} \subset m_+\mathcal{Q}_U^{-}$,
\item[$(ii)$] $\mathcal{P}^+ \cap \mathcal{M} = \{m_{+}\}$,
\item[$(iii)$] for all $m\in \mathcal{M}$ and all $(i,a)\in U$, if $mA_{i,a}^{-1} \not\in \mathcal{M}$, then
$mA_{i,a}^{-1}A_{j,b} \not\in \mathcal{M}$ unless $(j,b)=(i,a)$,
\item[$(iv)$] for all $m \in \mathcal{M}$ and all $i \in I$, there exists a~unique $i$-dominant monomial $M\in
\mathcal{M}$ such that
\begin{gather*}
\trunc_{\beta_i(M\mathcal{Q}_U^{-})} \chi_q(\beta_i(M))=\sum_{m'\in m\mathcal{Q}_{\{i\}
\times\mathbb{C}^{\times}}}\beta_i(m').
\end{gather*}
\end{enumerate}
Then
\begin{gather*}
\trunc_{m_{+}\mathcal{Q}_U^{-}} \chi_q(m_+)=\sum_{m\in\mathcal{M}}m.
\end{gather*}
\end{Theorem}

{\sloppy Here by $\chi_q(\beta_i(M))$ we mean the $q$-character of the irreducible
$U_{q_i}(\hat{\mathfrak{sl}_2})$-module with \mbox{highest} weight monomial $\beta_i(M)$ and by
$\trunc_{\beta_i}(M\mathcal{Q}_{U}^{-})$ we mean keeping only the monomials of $\chi_q(\beta_i(M))$
in the set $\beta_i(M\mathcal{Q}_U^-)$.

}

\section{Main results}
\label{main results}

\subsection{First examples}

Without loss of generality, we f\/ix $a\in \mathbb{C}^{\times}$ and consider modules $V$ with
$\mathscr{M}(V) \subset \mathbb{Z}[Y_{i, aq^{k}}]_{i\in I,\; k\in \mathbb{Z}}$.
In the following, for simplicity we write $i_s$, $i_s^{-1}$ ($s\in \mathbb{Z}$) instead of $Y_{i, aq^s}$,
$Y_{i, aq^s}^{-1}$ respectively.
The $q$-characters of fundamental modules are easy to compute by using the FM algorithm.
\begin{Lemma}
\label{fundamental modules}
The fundamental $q$-characters for $U_q\hat{\mathfrak{g}}$ of type $G_2$ are given by
\begin{gather*}
\chi_q(1_0)=1_0+1_2^{-1}2_1+1_41_62_7^{-1}+1_41_8^{-1}+1_6^{-1}1_8^{-1}2_5+1_{10}2_{11}^{-1}+1_{12}^{-1},
\\
\chi_q(2_0)=2_{0}+1_{1}1_{3}1_{5}2^{-1}_{6}+1_{1}1_{3}1^{-1}_{7}+1_{1}1^{-1}_{5}1^{-1}_{7}2_{4}+1^{-1}_{3}
1^{-1}_{5}1^{-1}_{7}2_{2}2_{4}
\\
\phantom{\chi_q(2_0)=}
{} + 1_{1}1_{9}2^{-1}_{10}+2_{4}2^{-1}_{8}+1^{-1}_{3}1_{9}2_{2}2^{-1}_{10}+1_{5}1_{7}1_{9}2^{-1}_{8}2^{-1}
_{10}+1_{1}1^{-1}_{11}
\\
\phantom{\chi_q(2_0)=}
{} + 1^{-1}_{3}1^{-1}_{11}2_{2}+1_{5}1_{7}1^{-1}_{11}2^{-1}_{8}+1_{5}1^{-1}_{9}1^{-1}_{11}+1^{-1}_{7}1^{-1}
_{9}1^{-1}_{11}2_{6}+2^{-1}_{12}.
\end{gather*}
\end{Lemma}

For $s\in \mathbb{Z}$, $\chi_q(1_s)$ and $\chi_q(2_s)$ are obtained by shift all indices by $s$ in
$\chi_q(1_0)$ and $\chi_q(2_0)$ respectively.

It is convenient to keep in mind the following lemma.
\begin{Lemma}
If $b\in \mathbb{Z} \backslash \{\pm 2, \pm 8, \pm 12\}$, then
\begin{gather*}
L(1_01_b)=L(1_0)\otimes L(1_b),
\qquad
\dim L(1_01_b)=49.
\end{gather*}
If $b\in \mathbb{Z} \backslash \{\pm 6, \pm 8, \pm 10, \pm 12\}$, then
\begin{gather*}
L(2_02_b)=L(2_0)\otimes L(2_b),
\qquad
\dim L(2_02_b)=225.
\end{gather*}
If $b\in \mathbb{Z} \backslash \{\pm 7, \pm 11\}$, then
\begin{gather*}
L(1_02_b)=L(1_0)\otimes L(2_b),
\qquad
L(2_01_b)=L(2_0)\otimes L(1_b),
\\
\dim L(1_02_b)=\dim L(2_01_b)=105.
\end{gather*}
In addition, we have
\begin{gather*}
\dim L(1_01_2)=34,
\qquad
\dim L(1_01_8)=42,
\qquad
\dim L(1_01_{12})=48,
\\
\dim L(2_02_6)=92,
\!\!\qquad
\dim L(2_02_8)=210,
\!\!\qquad
\dim L(2_02_{10})=183,
\!\!\qquad
\dim L(2_02_{12})=224,
\\
\dim L(1_02_7)=\dim L(2_01_7)=71,
\qquad
\dim L(1_02_{11})=\dim L(2_01_{11})=98.
\end{gather*}
\end{Lemma}

\begin{proof}
By Lemma~\ref{fundamental modules}, the tensor products in the f\/irst three cases of the lemma are special.
Therefore the tensor products are irreducible.
Hence the f\/irst three cases of the lemma are true.
The last part of the lemma can be proved using the methods of Section~\ref{proof system}.
In fact some of the dimensions follow from Theorem~\ref{dimensions}.
We do not use this lemma in the proofs.
Therefore we omit the details of the proof.
\end{proof}

\subsection[Definition of the modules $\mathcal{B}_{k, \ell}^{(s)}$, $\mathcal{C}_{k, \ell}^{(s)}$,
$\mathcal{D}_{k, \ell}^{(s)}$, $\mathcal{E}_{k, \ell}^{(s)}$, $\mathcal{F}_{k, \ell}^{(s)}$]
{Definition of the modules $\boldsymbol{\mathcal{B}_{k, \ell}^{(s)}}$, $\boldsymbol{\mathcal{C}_{k, \ell}^{(s)}}$,
$\boldsymbol{\mathcal{D}_{k, \ell}^{(s)}}$, $\boldsymbol{\mathcal{E}_{k, \ell}^{(s)}}$,
$\boldsymbol{\mathcal{F}_{k, \ell}^{(s)}}$}

For $s\in \mathbb{Z}$, $k,\ell\in \mathbb{Z}_{\geq 0}$, def\/ine the following monomials.
\begin{gather*}
B_{k,\ell}^{(s)}=\left(\prod_{i=0}^{k-1}2_{s+6i}\right)\left(\prod_{i=0}^{\ell-1}1_{s+6k+2i+1}
\right),
\qquad
C_{k,\ell}^{(s)}=\left(\prod_{i=0}^{k-1}2_{s+6i}\right)\left(\prod_{i=0}^{\ell-1}2_{s+6k+6i+4}
\right),
\\
D_{k,\ell}^{(s)}=\left(\prod_{i=0}^{k-1}2_{s+6i}\right)1_{s+6k+1}\left(\prod_{i=0}^{\ell-1}2_{s+6k+6i+8}
\right),
\qquad
F_{k,\ell}^{(s)}=\left(\prod_{i=0}^{k-1}1_{s+2i}\right)\!\left(\prod_{i=0}^{\ell-1}1_{s+2k+2i+6}
\right),
\\
E_{k,\ell}^{(s)}=\left(\prod_{i=0}^{k-1}1_{s+2i}\right)\left(\prod_{i=0}^{\left\lfloor\frac{\ell-1}{2}\right\rfloor}
2_{s+2k+6i+3}\right)\left(\prod_{i=0}^{\left\lfloor\frac{\ell-2}{2}\right\rfloor}2_{s+2k+6i+5}\right).
\end{gather*}

Note that, in particular, for $k \in \mathbb{Z}_{\geq 0}$, $s\in \mathbb{Z}$, we have the following trivial
relations
\begin{gather}
\mathcal{B}_{k,0}^{(s)}=\mathcal{C}_{k,0}^{(s)}=\mathcal{C}_{0,k}^{(s-4)},
\qquad
\mathcal{D}_{k,0}^{(s)}
=\mathcal{B}_{k,1}^{(s)},
\qquad
\mathcal{E}_{k,0}^{(s)}=\mathcal{B}_{0,k}^{(s-1)}=\mathcal{F}_{0,k}^{(s-6)}
=\mathcal{F}_{k,0}^{(s)}.
\label{trivial relations 1}
\end{gather}

Denote by $\mathcal{B}_{k, \ell}^{(s)}$, $\mathcal{C}_{k, \ell}^{(s)}$, $\mathcal{D}_{k, \ell}^{(s)}$,
$\mathcal{E}_{k, \ell}^{(s)}$, $\mathcal{F}_{k, \ell}^{(s)}$ the irreducible f\/inite-dimensional highest
$l$-weight $U_q\hat{\mathfrak{g}}$-modules with highest $l$-weight $B_{k, \ell}^{(s)}, C_{k, \ell}^{(s)},
D_{k, \ell}^{(s)}, E_{k, \ell}^{(s)}, F_{k, \ell}^{(s)}$ respectively.

Note that $\mathcal{B}_{k, \ell}^{(s)}$, $\mathcal{D}_{0, \ell}^{(s)}$, $\mathcal{D}_{k, 0}^{(s)}$ are
minimal af\/f\/inizations.
The modules $\mathcal{B}_{0, \ell}^{(s)}$, $\mathcal{C}_{0, \ell}^{(s)}$, $\mathcal{F}_{0, \ell}^{(s)}$,
$\mathcal{B}_{k, 0}^{(s)}$, $\mathcal{C}_{k, 0}^{(s)}$, $\mathcal{E}_{k, 0}^{(s)}$, $\mathcal{F}_{k,
0}^{(s)}$ are Kirillov--Reshetikhin modules.

Our f\/irst result is
\begin{Theorem}
\label{special}
The modules $\mathcal{B}_{k, \ell}^{(s)}$, $\mathcal{C}_{k, \ell}^{(s)}$, $\mathcal{D}_{k, \ell}^{(s)}$,
$\mathcal{E}_{k, \ell}^{(s)}$, $\mathcal{F}_{k, \ell}^{(s)}$, $s\in \mathbb{Z}$, $k,\ell\in \mathbb{Z}_{\geq0}$,
are special.
In particular, the FM algorithm works for all these modules.
\end{Theorem}

We prove Theorem~\ref{special} in Section~\ref{proof special}.
Note that the case of $\mathcal{B}_{k, \ell}^{(s)}$ has been proved in Theorem~3.8 of~\cite{Her07}.
In general, the modules in Theorem~\ref{special} are not thin.
For example, $\chi_q(1_01_2)$ has monomial $1_{4}1_{6}1_{8}^{-1}1_{10}^{-1}$ with coef\/f\/icient~$2$.

\subsection[Extended $T$-system]{Extended $\boldsymbol{T}$-system}

It is known that Kirillov--Reshetikhin modules $\mathcal{B}_{k, 0}^{(s)}$, $\mathcal{B}_{0, \ell}^{(s)}$
satisfy the following $T$-system relations, see~\cite{KNS94},
\begin{gather}
\Big[\mathcal{B}_{0,\ell}^{(s)}\Big]\Big[\mathcal{B}_{0,\ell}^{(s+2)}\Big]
=\Big[\mathcal{B}_{0,\ell+1}^{(s)}\Big]\Big[\mathcal{B}_{0,\ell-1}^{(s+2)}\Big]
+\Big[\mathcal{B}_{\left\lfloor\frac{\ell+2}{3}\right\rfloor,0}^{(s+1)}\Big]
\Big[\mathcal{B}_{\left\lfloor\frac{\ell+1}{3}\right\rfloor,0}^{(s+3)}\Big]
\Big[\mathcal{B}_{\left\lfloor\frac{\ell}{3}\right\rfloor,0}^{(s+5)}\Big],
\label{T-sys 1}
\\
\Big[\mathcal{B}_{k,0}^{(s)}\Big]\Big[\mathcal{B}_{k,0}^{(s+6)}\Big]
=\Big[\mathcal{B}_{k+1,0}^{(s)}\Big]\Big[\mathcal{B}_{k-1,0}
^{(s+6)}\Big]+\Big[\mathcal{B}_{0,3k}^{(s+1)}\Big],
\label{T-sys 2}
\end{gather}
where $s\in \mathbb{Z}$, $k,\ell\in \mathbb{Z}_{\geq 1}$.

Our main result is
\begin{Theorem}
\label{extended T-system}
For $s\in \mathbb{Z}$ and $k,\ell \in\mathbb{Z}_{\geq 1}$, $t\in\mathbb{Z}_{\geq 2}$, we have the following
relations in $\rep(U_q\hat{\mathfrak{g}})$:
\begin{gather}
\Big[\mathcal{B}_{k,\ell-1}^{(s)}\Big]\Big[\mathcal{B}_{k-1,\ell}^{(s+6)}\Big]
=\Big[\mathcal{B}_{k,\ell}^{(s)}\Big]\Big[\mathcal{B}_{k-1,\ell-1}^{(s+6)}\Big]
+\Big[\mathcal{E}_{3k-1,\left\lceil\frac{2\ell-2}{3}\right\rceil}^{(s+1)}\Big]
\Big[\mathcal{B}_{\left\lfloor\frac{\ell-1}{3}\right\rfloor,0}^{(s+6k+6)}\Big],
\label{b equation}
\\
\Big[\mathcal{E}_{0,\ell}^{(s)}\Big]
=\Big[\mathcal{B}_{\left\lfloor\frac{\ell+1}{2}\right\rfloor,0}^{(s+3)}\Big]
\Big[\mathcal{B}_{\left\lfloor\frac{\ell}{2}\right\rfloor,0}^{(s+5)}\Big],
\label{e equation 1}
\\
\Big[\mathcal{E}_{1, \ell}^{(s)}\Big]
=\Big[\mathcal{D}_{0,\left\lfloor\frac{\ell}{2}\right\rfloor}^{(s-1)}\Big]
\Big[\mathcal{B}_{\left\lfloor\frac{\ell+1}{2}\right\rfloor,0}^{(s+5)}\Big],
\label{e equation 2}
\\
\Big[\mathcal{E}_{t,\ell-1}^{(s)}\Big]\Big[\mathcal{E}_{t-1,\ell}^{(s+2)}\Big]
=\Big[\mathcal{E}_{t,\ell}^{(s)}\Big]\Big[\mathcal{E}
_{t-1,\ell-1}^{(s+2)}\Big]
\nonumber
\\
\hphantom{\Big[\mathcal{E}_{t,\ell-1}^{(s)}\Big]\Big[\mathcal{E}_{t-1,\ell}^{(s+2)}\Big]=}{}
+
\begin{cases}
\Big[\mathcal{D}_{r,p-1}^{(s+1)}\Big]\Big[\mathcal{B}_{r+p,0}^{(s+3)}\Big]\Big[\mathcal{B}_{r,3p-2}^{(s+5)}\Big],
&\text{if}\quad
t=3r+2,\; \ell=2p-1,
\\[2mm]
\Big[\mathcal{B}_{r+p+1,0}^{(s+1)}\Big]\Big[\mathcal{C}_{r,p}^{(s+3)}\Big]\Big[\mathcal{B}_{r,3p-1}^{(s+5)}\Big],
&\text{if}\quad
t=3r+2,\; \ell=2p,
\\[2mm]
\Big[\mathcal{B}_{r+1,3p-2}^{(s+1)}\Big]\Big[\mathcal{D}_{r,p-1}^{(s+3)}\Big]\Big[\mathcal{B}_{r+p,0}^{(s+5)}\Big],
&\text{if}\quad
t=3r+3,\; \ell=2p-1,
\\[2mm]
\Big[\mathcal{B}_{r+1,3p-1}^{(s+1)}\Big]\Big[\mathcal{B}_{r+p+1,0}^{(s+3)}\Big]\Big[\mathcal{C}_{r,p}^{(s+5)}\Big],
&\text{if}\quad
t=3r+3,\; \ell=2p,
\\[2mm]
\Big[\mathcal{B}_{r+p+1,0}^{(s+1)}\Big]\Big[\mathcal{B}_{r+1,3p-2}^{(s+3)}\Big]\Big[\mathcal{D}_{r,p-1}^{(s+5)}\Big],
&\text{if}\quad
t=3r+4,\; \ell=2p-1,
\\[2mm]
\Big[\mathcal{C}_{r+1,p}^{(s+1)}\Big]\Big[\mathcal{B}_{r+1,3p-1}^{(s+3)}\Big]\Big[\mathcal{B}_{r+p+1,0}^{(s+5)}\Big],
&\text{if}\quad
t=3r+4,\; \ell=2p,
\end{cases}
\label{e equation}
\\
\label{c equation}
\Big[\mathcal{C}_{k,\ell-1}^{(s)}\Big]\Big[\mathcal{C}_{k-1,\ell}^{(s+6)}\Big]
=\Big[\mathcal{C}_{k,\ell}^{(s)}\Big]\Big[\mathcal{C}
_{k-1,\ell-1}^{(s+6)}\Big]+\Big[\mathcal{F}_{3k-2,3\ell-2}^{(s+1)}\Big],
\\
\Big[\mathcal{D}_{0,\ell-1}^{(s)}\Big]\Big[\mathcal{B}_{\ell,0}^{(s+8)}\Big]
=\Big[\mathcal{D}_{0,\ell}^{(s)}\Big]\Big[\mathcal{B}
_{\ell-1,0}^{(s+8)}\Big]+\Big[\mathcal{B}_{0,3\ell-1}^{(s+4)}\Big],
\label{d equation 1}
\\
\label{d equation 2}
\Big[\mathcal{D}_{k,\ell-1}^{(s)}\Big]\Big[\mathcal{D}_{k-1,\ell}^{(s+6)}\Big]
=\Big[\mathcal{D}_{k,\ell}^{(s)}\Big]\Big[\mathcal{D}
_{k-1,\ell-1}^{(s+6)}\Big]+\Big[\mathcal{F}_{3k-1,3\ell-1}^{(s+1)}\Big],
\\
\Big[\mathcal{F}_{k,\ell-1}^{(s)}\Big]\Big[\mathcal{F}_{k-1,\ell}^{(s+2)}\Big]
=\Big[\mathcal{F}_{k,\ell}^{(s)}\Big]\Big[\mathcal{F}
_{k-1,\ell-1}^{(s+2)}\Big]
\label{f equation}
\\
\hphantom{\Big[\mathcal{F}_{k,\ell-1}^{(s)}\Big]\Big[\mathcal{F}_{k-1,\ell}^{(s+2)}\Big]=}{}
+
\begin{cases}
\Big[\mathcal{B}_{r,0}^{(s+1)}\Big]\Big[\mathcal{D}_{r,\left\lfloor\frac{\ell}{3}\right\rfloor}^{(s+3)}\Big]
\Big[\mathcal{C}_{r,\left\lfloor\frac{\ell+1}{3}\right\rfloor}^{(s+5)}\Big]
\Big[\mathcal{B}_{\left\lfloor\frac{\ell-1}{3}\right\rfloor,0}^{(s+2k+11)}\Big],&\text{if}\quad k=3r+1,
\\[2mm]
\Big[\mathcal{C}_{r+1,\left\lfloor\frac{\ell+1}{3}\right\rfloor}^{(s+1)}\Big]\Big[\mathcal{B}_{r,0}^{(s+3)}\Big]
\Big[\mathcal{D}_{r,\left\lfloor\frac{\ell}{3}\right\rfloor}^{(s+5)}\Big]
\Big[\mathcal{B}_{\left\lfloor\frac{\ell-1}{3}\right\rfloor,0}^{(s+2k+11)}\Big],&\text{if}\quad k=3r+2,
\\[2mm]
\Big[\mathcal{D}_{r+1,\left\lfloor\frac{\ell}{3}\right\rfloor}^{(s+1)}\Big]
\Big[\mathcal{C}_{r+1,\left\lfloor\frac{\ell+1}{3}\right\rfloor}^{(s+3)}\Big]
\Big[\mathcal{B}_{r,0}^{(s+5)}\Big]\Big[\mathcal{B}_{\left\lfloor\frac{\ell-1}{3}\right\rfloor,0}^{(s+2k+11)}\Big],
&\text{if}\quad
k=3r+3.
\end{cases} \nonumber
\end{gather}
\end{Theorem}

We prove Theorem~\ref{extended T-system} in Section~\ref{proof system}.

Note that since $D_{k, 0}^{(s)} = B_{k, 1}^{(s)}$, equations for $\mathcal{D}_{k, 0}^{(s)}$ are included in
the equations for $\mathcal{B}_{k, 1}^{(s)}$.

All relations except~\eqref{e equation 1},~\eqref{e equation 2} in Theorem~\ref{extended T-system} are
written in the form $[\mathcal{L}] [\mathcal{R}]=[\mathcal{T}] [\mathcal{B}] + [\mathcal{S}]$, where
$\mathcal{L}$, $\mathcal{R}$, $\mathcal{T}$, $\mathcal{B}$ are irreducible modules which we call \textit{left,
right, top and bottom modules} and~$\mathcal{S}$ is a~tensor product of some irreducible modules.
We call the factors of~$\mathcal{S}$ \textit{sources}.
Moreover, we have the following theorem.
\begin{Theorem}
\label{irreducible}
For each relation in Theorem~{\rm \ref{extended T-system}}, all summands on the right hand
side, $\mathcal{T}\otimes \mathcal{B}$ and $\mathcal{S}$, are irreducible.
\end{Theorem}

We will prove Theorem~\ref{irreducible} in Section~\ref{proof irreducible}.

Recall that the $q$-characters of modules for dif\/ferent $s$ are related by the simple shift of indexes,
see~\eqref{shift}.

We have the following proposition.
\begin{Proposition}
\label{compute}
Given $\chi_q(1_s)$, $\chi_q(2_s)$, one can obtain the $q$-characters of $\mathcal{B}_{k, \ell}^{(s)}$,
$\mathcal{C}_{k, \ell}^{(s)}$, $\mathcal{D}_{k, \ell}^{(s)}$, $\mathcal{E}_{k, \ell}^{(s)}$,
$\mathcal{F}_{k, \ell}^{(s)}$, $s\in\mathbb{Z}$, $k,\ell\in \mathbb{Z}_{\geq 0}$, recursively, by
using~\eqref{trivial relations 1}, and computing the $q$-character of the top module through the
$q$-characters of other modules in relations in Theorem~{\rm \ref{extended T-system}}.
\end{Proposition}
\begin{proof}
\textbf{Claim 1.} Let $n$, $m$ be positive integers.
Then the $q$-characters
\begin{gather*}
\chi_q\Big(\mathcal{B}_{k,\ell}^{(s)}\Big),
\quad
k\leq n,
\quad
\ell\leq m,
\qquad
\chi_q\Big(\mathcal{C}_{k,\ell}^{(s)}\Big),
\quad
k\leq n-1,
\quad
\ell \leq\left\lceil\frac{2m+1}{6}\right\rceil,
\\
\chi_q\Big(\mathcal{D}_{k,\ell}^{(s)}\Big),
\quad
k\leq n-1,
\quad
\ell\leq\left\lceil\frac{2m+1}{6}\right\rceil,
\qquad
\chi_q\Big(\mathcal{E}_{k,\ell}^{(s)}\Big),
\quad
k\leq3n-1,
\quad
\ell\leq\left\lceil\frac{2m-2}{3}\right\rceil,
\\
\chi_q\Big(\mathcal{F}_{k,\ell}^{(s)}\Big),
\quad
k\leq3n-4,
\quad
\ell\leq m+2,
\end{gather*}
can be computed recursively starting from $\chi_q(1_0)$, $\chi_q(2_0)$.

We use induction on $n$, $m$ to prove Claim~1.
For simplicity, we do not write the upper-subscripts ``$(s)$'' in the remaining part of the proof.
We know that, see~\cite{Her06}, the $q$-characters of Kirillov--Reshetikhin modules can be computed from
$\chi_q(1_0)$, $\chi_q(2_0)$.

When $n=0$, $m=1$, Claim~1 is clearly true.
It is clear that $\chi_q(\mathcal{D}_{0,1})$ can be computed using~\eqref{d equation 1}.
Therefore Claim~1 holds for $n=1$, $m=0$,

Suppose that for $n\leq n_1$ and $m \leq m_1$, Claim~1 is true.
Let $n=n_1+1$, $m=m_1$.
We need to show that Claim~1 is true.
Then we need to show that
\begin{gather*}
\chi_q(\mathcal{B}_{n_1+1,\ell}),
\quad
\ell\leq m_1,
\qquad
\chi_q(\mathcal{C}_{n_1,\ell}),
\quad
\ell\leq\left\lceil\frac{2m_1+1}{6}\right\rceil,
\qquad
\chi_q(\mathcal{D}_{n_1,\ell}),
\quad
\ell\leq\left\lceil\frac{2m_1+1}{6}\right\rceil,
\\
\chi_q(\mathcal{E}_{k,\ell}),
\quad
k=3n_1,
\quad
3n_1+1,
\quad
3n_1+2,
\quad
\ell\leq\left\lceil\frac{2m_1-2}{3}\right\rceil,
\\
\chi_q(\mathcal{F}_{k,\ell}),
\quad
k=3n_1-3,
\quad
3n_1-2,
\quad
3n_1-1,
\quad \ell\leq m_1+2,
\end{gather*}
can be computed.

We compute the following modules
\begin{gather*}
\chi_q(\mathcal{F}_{3n_1-3,\ell}),
\quad
\ell\leq m_1+2,
\qquad
\chi_q(\mathcal{F}_{3n_1-2,\ell}),
\quad
\ell\leq m_1+2,
\\
\chi_q(\mathcal{C}_{n_1,\ell}),
\quad
\ell\leq\left\lfloor\frac{m_1+3}{3}\right\rfloor,
\qquad
\chi_q(\mathcal{F}_{3n_1-1,\ell}),
\quad
\ell\leq m_1+2,
\\
\chi_q(\mathcal{D}_{n_1,\ell}),
\quad
\ell\leq\left\lceil\frac{2m_1+1}{6}\right\rceil,
\qquad
\chi_q(\mathcal{C}_{n_1,\ell}),
\quad
\ell\leq\left\lceil\frac{2m_1+1}{6}\right\rceil,
\\
\chi_q(\mathcal{E}_{3n_1,\ell}),
\quad
\ell\leq\left\lceil\frac{2m_1-2}{3}\right\rceil,
\qquad
\chi_q(\mathcal{E}_{3n_1,\ell}),
\quad
\ell\leq\left\lceil\frac{2m_1-2}{3}\right\rceil,
\\
\chi_q(\mathcal{E}_{3n_1+1,\ell}),
\quad
\ell\leq\left\lceil\frac{2m_1-2}{3}\right\rceil,
\qquad
\chi_q(\mathcal{E}_{3n_1+2,\ell}),
\quad
\ell\leq\left\lceil\frac{2m_1-2}{3}\right\rceil,
\\
\chi_q(\mathcal{B}_{n_1+1,\ell}),
\quad
\ell\leq m_1
\end{gather*}
in the order as shown.
At each step, we consider the module that we want to compute as a~top module and use the corresponding
relation in Theorem~\ref{extended T-system} and known $q$-characters.
For example, we consider the f\/irst set of modules $\chi_q(\mathcal{F}_{3n_1-3, \ell})$, $\ell \leq m_1+2$.
Since $\left\lfloor \frac{m_1+3}{3} \right\rfloor \leq \left\lceil \frac{2m_1+1}{6} \right\rceil$,
$\chi_q(\mathcal{C}_{n_1-1, \ell})$, $\ell \leq \left\lfloor \frac{m_1+3}{3} \right\rfloor$, is known by induction
hypothesis.
Similarly, $\chi_q(\mathcal{D}_{n_1-1, \ell})$, $\ell \leq \left\lfloor \frac{m_1+2}{3} \right\rfloor$ is known.
Therefore $\chi_q(\mathcal{F}_{3n_1-3, \ell})$, $\ell \leq m_1+2$, is computed using the last equation
of~\eqref{f equation}.

Similarly, we show that Claim 1 holds for $n=n_1$, $m=m_1+1$.
Therefore Claim 1 is true for all $n \geq 1$, $m \geq 1$.
\end{proof}

\section{Proof of Theorem~\ref{special}}
\label{proof special}

In this section, we will show that the modules $\mathcal{B}_{k, \ell}^{(s)}$, $\mathcal{C}_{k,
\ell}^{(s)}$, $\mathcal{D}_{k, \ell}^{(s)}$, $\mathcal{E}_{k, \ell}^{(s)}$, $\mathcal{F}_{k, \ell}^{(s)}$
are special.

Since $\mathcal{B}_{0, \ell}^{(s)}$, $\mathcal{C}_{0, \ell}^{(s)}$, $\mathcal{F}_{0, \ell}^{(s)}$,
$\mathcal{B}_{k, 0}^{(s)}$, $\mathcal{C}_{k, 0}^{(s)}$, $\mathcal{E}_{k, 0}^{(s)}$, $\mathcal{F}_{k,
0}^{(s)}$ are Kirillov--Reshetikhin modules, they are special.

\subsection[The case of $\mathcal{C}_{k, \ell}^{(s)}$]{The case of $\boldsymbol{\mathcal{C}_{k, \ell}^{(s)}}$}

Let $m_+=C_{k, \ell}^{(s)}$ with $k,\ell\in\mathbb{Z}_{\geq 1}$.
Without loss of generality, we can assume that $s=6$.
Then
\begin{gather*}
m_+=(2_{6}2_{12}\cdots2_{6k})(2_{6k+10}2_{6k+16}\cdots2_{6k+6\ell+4}).
\end{gather*}

\textbf{Case 1.} $k=1$.
Let $U=I \times \{aq^s: s \in \mathbb{Z}, s < 6\ell +13 \}$.
Clearly, all monomials in $\chi_q(m_+)-\trunc_{m_+ \mathcal{Q}_{U}^{-}}   \chi_q(m_+)$ are
right-negative.
Therefore it is suf\/f\/icient to show that $\trunc_{m_+ \mathcal{Q}_{U}^{-}}   \chi_q(m_+)$ is
special.

Let $\mathcal{M}$ be the f\/inite set consisting of the following monomials
\begin{gather*}
m_0=m_+,
\qquad
m_1=m_0A_{2,9}^{-1},
\qquad
m_2=m_1A_{1,12}^{-1},
\\
m_{3}=m_{2}A_{1,10}^{-1},
\qquad
m_4=m_{3}A_{1,8}^{-1},
\qquad
m_5=m_{4}A_{2,11}^{-1}.
\end{gather*}
It is clear that $\mathcal{M}$ satisf\/ies the conditions in Theorem~\ref{truncated}.
Therefore
\begin{gather*}
\trunc_{m_+\mathcal{Q}_{U}^{-}} \chi_q(m_+)=\sum_{m\in\mathcal{M}}m
\end{gather*}
and $\trunc_{m_+ \mathcal{Q}_{U}^{-}}   \chi_q(m_+)$ is special.

\textbf{Case 2.} $k>1$.
Since the conditions of Theorem~\ref{truncated} do not apply to this case, we use another technique to show
that $L(m_+)$ is special.
We embed $L(m_+)$ into two dif\/ferent tensor products.
In both tensor products, each factor is special.
Therefore we can use the FM algorithm to compute the $q$-characters of the factors.
We classify the dominant monomials in the f\/irst tensor product and show that the only dominant monomial
in the f\/irst tensor product which occurs in the second tensor product is $m_+$ which proves that $L(m_+)$
is special.

The f\/irst tensor product is $L(m'_1) \otimes L(m'_2)$, where
\begin{gather*}
m'_1=2_{6}2_{12}\cdots2_{6k},
\qquad
m'_2=2_{6k+10}2_{6k+16}\cdots2_{6k+6\ell+4}.
\end{gather*}

We use the FM algorithm to compute $\chi_q(m'_1), \chi_q(m'_2)$ and classify all dominant monomials in
$\chi_q(m'_1) \chi_q(m'_2)$.
Let $m=m_1m_2$ be a~dominant monomial, where $m_i \in \chi_q(m'_i)$, $i=1, 2$.
If $m_2 \neq m'_2$, then $m$ is a~right negative monomial therefore $m$ is not dominant.
Hence $m_2=m'_2$.

If $m_1 \neq m'_1$, then $m_1$ is right negative.
Since $m$ is dominant, each factor with a~negative power in $m_1$ needs to be canceled by a~factor in
$m'_2$.
All possible cancellations cancel $2_{6k+10}$ in $m'_2$.
We have $\mathcal{M}(L(m'_1)) \subset \mathcal{M}(\chi_q(2_{6}2_{12}\cdots 2_{6k-6}) \chi_q(2_{6k}))$.
Only monomials in~$\chi_q(2_{6k})$ can cancel $2_{6k+10}$.
These monomials are $1_{6k+1}1_{6k+9}2_{6k+10}^{-1}$, $1_{6k+3}^{-1}1_{6k+9}2_{6k+2}2_{6k+10}^{-1}$, and
$1_{6k+5}1_{6k+7}1_{6k+9}2_{6k+8}^{-1}2_{6k+10}^{-1}$.
Therefore $m_1$ is in one of the following polynomials
\begin{gather}
\chi_q(2_{6}2_{12}\cdots2_{6k-6})1_{6k+1}1_{6k+9}2_{6k+10}^{-1},
\label{set 1}
\\
\chi_q(2_{6}2_{12}\cdots2_{6k-6})1_{6k+3}^{-1}1_{6k+9}2_{6k+2}2_{6k+10}^{-1},
\label{set 2}
\\
\chi_q(2_{6}2_{12}\cdots2_{6k-6})1_{6k+5}1_{6k+7}1_{6k+9}2_{6k+8}^{-1}2_{6k+10}^{-1}.
\label{set 3}
\end{gather}

\textbf{Subcase 2.1.} Let $m_1$ be in~\eqref{set 1}.
If $m_1 = 2_{6}2_{12}\cdots 2_{6k-6} 1_{6k+1}1_{6k+9}2_{6k+10}^{-1}$, then
\begin{gather}
m=m_1m_2=2_{6}2_{12}\cdots2_{6k-6}1_{6k+1}1_{6k+9}2_{6k+16}\cdots 2_{6k+6\ell+4}
\label{the other dominant monomial}
\end{gather}
is dominant.
Suppose that
\begin{gather*}
m_1\neq2_{6}2_{12}\cdots2_{6k-6}1_{6k+1}1_{6k+9}2_{6k+10}^{-1}.
\end{gather*}
Then $m_1=n_1 1_{6k+1}1_{6k+9}2_{6k+10}^{-1}$, where $n_1$ is a~non-highest monomial in
$\chi_q(2_{6}2_{12}\cdots 2_{6k-6})$.
Since $n_1$ is right negative, $1_{6k+1}$ or $1_{6k+9}$ should cancel a~factor of $n_1$ with a~negative
power.
Using the FM algorithm, we see that there exists a~factor $1_{6k-1}^2$ or $1_{6k+7}^{2}$ in a~monomial in
$\chi_q(2_{6}2_{12}\cdots 2_{6k-6})1_{6k+1}1_{6k+9}2_{6k+10}^{-1}$.
By Lemma~\ref{fundamental modules}, neither $1_{6k-1}^2$ nor $1_{6k+7}^{2}$ appear.
This is a~contradiction.

\textbf{Subcase 2.2.} Let $m_1$ be in~\eqref{set 2}.
If $m_1 = 2_{6}2_{12}\cdots 2_{6k-6} 1_{6k+3}^{-1}1_{6k+9}2_{6k+2}2_{6k+10}^{-1}$, then $m=m_1m_2$ is not
dominant.
Suppose that $m_1 \neq 2_{6}2_{12}\cdots 2_{6k-6} 1_{6k+3}^{-1}1_{6k+9}2_{6k+2}2_{6k+10}^{-1}$.
Then $m_1=n_1 1_{6k+3}^{-1}1_{6k+9}2_{6k+2}2_{6k+10}^{-1}$, where $n_1$ is a~non-highest monomial in
$\chi_q(2_{6}2_{12}\cdots 2_{6k-6})$.
Since $n_1$ is right negative, $1_{6k+9}$ or $2_{6k+2}$ should cancel a~factor of $n_1$ with a~negative
power.
Using the FM algorithm, we see that there exists either a~factor $1_{6k+7}^{2}$ or a~factor $2_{6k-4}^2$ in
a~monomial in $\chi_q(2_{6}2_{12}\cdots 2_{6k-6})1_{6k+1}1_{6k+9}2_{6k+10}^{-1}$.
By Lemma~\ref{fundamental modules}, neither $1_{6k+7}^{2}$ nor $2_{6k-4}^2$ appear.
This is a~contradiction.

\textbf{Subcase 2.3.} Let $m_1$ be in~\eqref{set 3}.
If $m_1 = 2_{6}2_{12}\cdots 2_{6k-6} 1_{6k+5}1_{6k+7}1_{6k+9}2_{6k+8}^{-1}2_{6k+10}^{-1}$, then $m=m_1m_2$
is not dominant.
Suppose that $m_1 \neq 2_{6}2_{12}\cdots 2_{6k-6} 1_{6k+5}1_{6k+7}1_{6k+9}2_{6k+8}^{-1}2_{6k+10}^{-1}$.
Then we have $m_1=n_1 1_{6k+5}1_{6k+7}1_{6k+9}2_{6k+8}^{-1}2_{6k+10}^{-1}$, where $n_1$ is a~non-highest
monomial in $\chi_q(2_{6}2_{12}\cdots 2_{6k-6})$.
Since $n_1$ is right negative, $1_{6k+5}$ or $1_{6k+7}$ or $1_{6k+9}$ should cancel a~factor of $n_1$ with
a~negative power.
Using the FM algorithm, we see that there exists a~factor $1_{6k+7}$ or $1_{6k+5}$ or $1_{6k+3}^{2}$ in
a~monomial in $\chi_q(2_{6}2_{12}\cdots 2_{6k-6})1_{6k+1}1_{6k+9}2_{6k+10}^{-1}$.
By Lemma~\ref{fundamental modules}, $1_{6k+7}$, $1_{6k+5}$, and $1_{6k+3}^{2}$ do not appear.
This is a~contradiction.

Therefore the only dominant monomials in $\chi_q(m'_1) \chi_q(m'_2)$ are $m_+$ and~\eqref{the other
dominant monomial}.

The second tensor product is $L(m''_1) \otimes L(m''_2)$, where
\begin{gather*}
m''_1=2_{6}2_{12}\cdots 2_{6k-6},
\qquad
m''_2=2_{6k}2_{6k+10}2_{6k+16}\cdots 2_{6k+6\ell+4}.
\end{gather*}

The monomial~\eqref{the other dominant monomial} is
\begin{gather}
n=m_+A_{2,6k+3}^{-1}A_{1,6k+6}^{-1}A_{1,6k+4}^{-1}A_{2,6k+7}^{-1}.
\label{expression of n}
\end{gather}
Since $A_{i, a}$, $i\in I$, $a\in \mathbb{C}^{\times}$ are algebraically independent, the
expression~\eqref{expression of n} of $n$ of the form $m_+\prod_{i\in I,\;a\in\mathbb{C}^{\times}}
A_{i, a}^{-v_{i, a}}$, where $v_{i, a}$ are some integers, is unique.
Suppose that the monomial $n$ is in $\chi_q(m''_1) \chi_q(m''_2)$.
Then $n=n_1n_2$, where $n_i \in \chi_q(m''_i)$, $i=1, 2$.
By the expression~\eqref{expression of n}, we have $n_1=m''_1$ and
\begin{gather*}
n_2=m''_2A_{2,6k+3}^{-1}A_{1,6k+6}^{-1}A_{1,6k+4}^{-1}A_{2,6k+7}^{-1}.
\end{gather*}
By the FM algorithm, the monomial $m''_2A_{2, 6k+3}^{-1}A_{1, 6k+6}^{-1}A_{1, 6k+4}^{-1}A_{2, 6k+7}^{-1}$
is not in $\chi_q(m''_2)$.
This contradicts the fact that $n_2 \in \chi_q(m''_2)$.
Therefore $n$ is not in $\chi_q(m''_1) \chi_q(m''_2)$.

\subsection[The case of $\mathcal{B}_{k, \ell}^{(s)}$]{The case of $\boldsymbol{\mathcal{B}_{k, \ell}^{(s)}}$}

Let $m_+=B_{k, \ell}^{(s)}$ with $k$, $\ell\in\mathbb{Z}_{\geq 1}$.
Without loss of generality, we can assume that $s=6$.
Then
\begin{gather*}
m_+=(2_{6}2_{12}\cdots2_{6k})(1_{6k+7}1_{6k+9}\cdots1_{6k+2\ell+5}).
\end{gather*}
Let $U=I \times \{aq^s: s\in \mathbb{Z}, s < 6k+2\ell +6\}$.
Clearly, all monomials in the polynomial $\chi_q(m_+)-\trunc_{m_+ \mathcal{Q}_{U}^{-}}
\chi_q(m_+)$ are right-negative.
Therefore it is suf\/f\/icient to show that the truncated $q$-character $\trunc_{m_+
\mathcal{Q}_{U}^{-}}   \chi_q(m_+)$ is special.

Let $\mathcal{M}$ be the f\/inite set consisting of the following monomials
\begin{gather*}
m_0=m_+,
\qquad
m_1=m_0A_{2,6k+3}^{-1},
\quad
m_2=m_1A_{2,6k-3}^{-1},
\quad
\ldots,
\quad
m_{k}=m_{k-1}A_{2,9}^{-1}.
\end{gather*}
It is clear that $\mathcal{M}$ satisf\/ies the conditions in Theorem~\ref{truncated}.
Therefore
\begin{gather*}
\trunc_{m_+\mathcal{Q}_{U}^{-}} \chi_q(m_+)=\sum_{m\in\mathcal{M}}m
\end{gather*}
and $\trunc_{m_+ \mathcal{Q}_{U}^{-}}   \chi_q(m_+)$ is special.

\subsection[The case of $\mathcal{D}_{k, \ell}^{(s)}$]{The case of $\boldsymbol{\mathcal{D}_{k, \ell}^{(s)}}$}

Let $m_+=D_{k, \ell}^{(s)}$ with $k$, $\ell \in\mathbb{Z}_{\geq 0}$.
Without loss of generality, we can assume that $s=0$.
Then
\begin{gather*}
m_+=(2_{0}2_{6}\cdots2_{6k-6})1_{6k+1}(2_{6k+8}2_{6k+14}\cdots2_{6k+6\ell+2}).
\end{gather*}

\textbf{Case 1.} $k = 0$.
Let $U=I \times \{aq^s: s\in \mathbb{Z}, s < 6\ell + 5 \}$.
Clearly, all monomials in $\chi_q(m_+)-\trunc_{m_+ \mathcal{Q}_{U}^{-}}   \chi_q(m_+)$ are
right-negative.
Therefore it is suf\/f\/icient to show that $\trunc_{m_+ \mathcal{Q}_{U}^{-}}   \chi_q(m_+)$ is
special.

Let
\begin{gather*}
M=\big\{m_+,m_{+}A_{1,2}^{-1}\big\}.
\end{gather*}
It is clear that $\mathcal{M}$ satisf\/ies the conditions in Theorem~\ref{truncated}.
Therefore
\begin{gather*}
\trunc_{m_+\mathcal{Q}_{U}^{-}} \chi_q(m_+)=\sum_{m\in\mathcal{M}}m
\end{gather*}
and $\trunc_{m_+ \mathcal{Q}_{U}^{-}}   \chi_q(m_+)$ is special.

\textbf{Case 2.} $k > 0$.
Let
\begin{gather*}
m'_1=2_{0}2_{6}\cdots2_{6k-6}1_{6k+1},
\qquad
m'_2=2_{6k+8}2_{6k+14}\cdots2_{6k+6\ell+2},
\\
m''_1=2_{0}2_{6}\cdots2_{6k-6},
\qquad
m''_2=1_{6k+1}2_{6k+8}2_{6k+14}\cdots2_{6k+6\ell+2}.
\end{gather*}
Then $\mathscr{M}(L(m_+)) \subset \mathscr{M}(\chi_q(m'_1)\chi_q(m'_2)) \cap
\mathscr{M}(\chi_q(m''_1)\chi_q(m''_2))$.

By using similar arguments as the case of $\mathcal{C}_{k, \ell}^{(s)}$, we show that the only dominant
monomials in $\chi_q(m'_1)\chi_q(m'_2)$ are $m_+$ and
\begin{gather*}
n=2_{0}2_{6}\cdots2_{6k-6}1_{6k+5}1_{6k+7}2_{6k+14}2_{6k+20}\cdots2_{6k+6\ell+2}=m_+A_{1,6k+2}^{-1}
A_{2,6k+5}^{-1}.
\end{gather*}
Moreover, $n$ is not in $\chi_q(m''_1)\chi_q(m''_2)$.
Therefore the only dominant monomial in $\chi_q(m_+)$ is $m_+$.

\subsection[The case of $\mathcal{E}_{k, \ell}^{(s)}$]{The case of $\boldsymbol{\mathcal{E}_{k, \ell}^{(s)}}$}

Let $m_{+}=E_{k, \ell}^{(s)}$ with $k$, $\ell\in\mathbb{Z}_{\geq 0}$.
Without loss of generality, we can assume that $s=1$.
Suppose that $\ell=2r+1$, $r\geq 0$ and $k = 3p$, $p\geq 1$.
The cases of $\ell=2r$, $r\geq 1$, or $k=0$ or $k=3p+1$, $p\geq 0$ or $k=3p+2$, $p\geq 0$ are similar.

Then
\begin{gather*}
m=(1_{1}1_{3}\cdots 1_{6p-1})(2_{6p+4}2_{6p+10}\cdots 2_{6p+6r-2}2_{6p+6r+4})(2_{6p+6}2_{6p+12}\cdots 2_{6p+6r}
).
\end{gather*}
Let $U=I \times \{aq^s: s\in \mathbb{Z}, s < 6p+6r+3\}$.
Clearly, all monomials in the polynomial $\chi_q(m_+)-\trunc_{m_+ \mathcal{Q}_{U}^{-}}
\chi_q(m_+)$ are right-negative.
Therefore it is suf\/f\/icient to show that the truncated $q$-character $\trunc_{m_+
\mathcal{Q}_{U}^{-}}   \chi_q(m_+)$ is special.

Let $\mathcal{M}$ be the f\/inite set consisting of the following monomials
\begin{gather*}
m_0=m_+,
\quad
m_1=m_0A_{1,6p}^{-1},
\quad
m_2=m_1A_{1,6p-2}^{-1},
\quad
\ldots,
\quad
m_{3p}=m_{3p-1}A_{1,2}^{-1},
\\
m_{3p+1}=m_{3p}A_{2,6p-4}^{-1},
\quad
m_{3p+2}=m_{3p+1}A_{2,6p-10}^{-1},
\quad
\ldots,
\quad
m_{4p}=m_{4p-2}A_{2,6}^{-1}.
\end{gather*}
It is clear that $\mathcal{M}$ satisf\/ies the conditions in Theorem~\ref{truncated}.
Therefore
\begin{gather*}
\trunc_{m_+\mathcal{Q}_{U}^{-}} \chi_q(m_+)=\sum_{m\in\mathcal{M}}m
\end{gather*}
and $\trunc_{m_+ \mathcal{Q}_{U}^{-}}   \chi_q(m_+)$ is special.

\subsection[The case of $\mathcal{F}_{k, \ell}^{(s)}$]{The case of $\boldsymbol{\mathcal{F}_{k, \ell}^{(s)}}$}

Let $m_+=F_{k, \ell}^{(s)}$ with $k$, $\ell\in\mathbb{Z}_{\geq 1}$.
Without loss of generality, we can assume that $s=1$.
Then
\begin{gather*}
m_+=(1_{1}1_{3}\cdots1_{2k-1})(1_{2k+7}1_{2k+9}\cdots1_{2k+2\ell+5}).
\end{gather*}

\textbf{Case 1.} $k=1$.
Let $U=I \times \{aq^{s}: s\in \mathbb{Z}, s < 2\ell + 8 \}$.
Clearly, all monomials in $\chi_q(m_+)-\trunc_{m_+ \mathcal{Q}_{U}^{-}}  \chi_q(m_+)$ are
right-negative.
Therefore it is suf\/f\/icient to show that $\trunc_{m_+ \mathcal{Q}_{U}^{-}}   \chi_q(m_+)$ is
special.

Let $\mathcal{M}$ be the f\/inite set consisting of the following monomials
\begin{gather*}
m_0=m_+,
\qquad
m_1=m_0A_{1,2}^{-1},
\qquad
m_2=m_1A_{2,5}^{-1}.
\end{gather*}
It is clear that $\mathcal{M}$ satisf\/ies the conditions in Theorem~\ref{truncated}.
Therefore
\begin{gather*}
\trunc_{m_+\mathcal{Q}_{U}^{-}} \chi_q(m_+)=\sum_{m\in\mathcal{M}}m
\end{gather*}
and $\trunc_{m_+ \mathcal{Q}_{U}^{-}}   \chi_q(m_+)$ is special.

\textbf{Case 2.} $k>1$.
Let
\begin{gather*}
m'_1=1_{1}1_{3}\cdots1_{2k-1},
\qquad
m'_2=1_{2k+7}1_{2k+9}\cdots1_{2k+2\ell+5},
\\
m''_1=1_{1}1_{3}\cdots1_{2k-3},
\qquad
m''_2=1_{2k-1}1_{2k+7}1_{2k+9}\cdots1_{2k+2\ell+5}.
\end{gather*}
Then $\mathscr{M}(L(m_+)) \subset \mathscr{M}(\chi_q(m'_1)\chi_q(m'_2)) \cap
\mathscr{M}(\chi_q(m''_1)\chi_q(m''_2))$.

By using similar arguments as the case of $\mathcal{C}_{k, \ell}^{(s)}$, we can show that the only dominant
monomials in $\chi_q(m'_1)\chi_q(m'_2)$ are $m_+$ and
\begin{gather*}
n_1=1_{1}1_{3}\cdots1_{2k-3}1_{2k+3}1_{2k+9}1_{2k+11}\cdots1_{2k+2\ell+5}
=m_+A_{1,2k}^{-1}A_{2,2k+3}^{-1}A_{1,2k+6}^{-1},
\\
n_2=1_{1}1_{3}\cdots1_{2k-3}1_{2k+7}1_{2k+9}1_{2k+13}1_{2k+15}\cdots1_{2k+2\ell+5}
=n_1A_{1,2k+4}^{-1}A_{2,2k+7}^{-1}A_{1,2k+10}^{-1},
\\
n_3=1_{1}1_{3}\cdots1_{2k-5}1_{2k+7}1_{2k+13}1_{2k+15}\cdots1_{2k+2\ell+5}
\\
\phantom{n_3}
=n_2A_{1,2k-2}^{-1}A_{2,2k+1}^{-1}A_{1,2k+4}^{-1}A_{1,2k+2}^{-1}A_{2,2k+5}^{-1}A_{1,2k+8}^{-1},
\\
n_4=1_{1}1_{3}\cdots1_{2k-7}1_{2k+13}1_{2k+15}\cdots1_{2k+2\ell+5}
\\
\phantom{n_4}
=n_3A_{1,2k-4}^{-1}A_{2,2k-1}^{-1}A_{1,2k+2}^{-1}A_{1,2k}^{-1}A_{2,2k+3}^{-1}A_{1,2k+6}^{-1}.
\end{gather*}
Moreover, $n_1$, $n_2$, $n_3$, $n_4$ are not in $\chi_q(m'_1)\chi_q(m'_2)$.
Therefore the only dominant monomial in~$\chi_q(m_+)$ is $m_+$.

\section{Proof of Theorem~\ref{extended T-system}}\label{proof system}

{\sloppy We use the FM algorithm to classify dominant monomials in $\chi_q(\mathcal{L})\chi_q(\mathcal{R})$,
$\chi_q(\mathcal{T})\chi_q(\mathcal{B})$, and~$\chi_q(\mathcal{S})$.

}

\subsection[Classif\/ication of dominant monomials in $\chi_q(\mathcal{L})\chi_q(\mathcal{R})$ and
$\chi_q(\mathcal{T})\chi_q(\mathcal{B})$]
{Classif\/ication of dominant monomials in $\boldsymbol{\chi_q(\mathcal{L})\chi_q(\mathcal{R})}$ and
$\boldsymbol{\chi_q(\mathcal{T})\chi_q(\mathcal{B})}$}

\begin{Lemma}
\label{dominant monomials}
We have the following cases.
\begin{enumerate}\itemsep=0pt
\item[$(1)$] Let $M=B_{k, \ell-1}^{(s)}B_{k-1, \ell}^{(s+6)}$, $k\geq 1$, $\ell \geq 1$.
Then dominant monomials
in $\chi_q\Big(\mathcal{B}_{k, \ell-1}^{(s)}\Big)\chi_q\Big(\mathcal{B}_{k-1,\ell}^{(s+6)}\Big)$ are
\begin{gather*}
M_0=M,
\quad
M_1=MA_{1,s+6k+2\ell-2}^{-1},
\\
M_2=M_1A_{1,s+6k+2\ell-4}^{-1},
\quad
\ldots,
\quad
M_{\ell-1}=M_{\ell-2}A_{1,s+6k+2}^{-1},
\\
M_{\ell}=M_{\ell-1}A_{2,s+6k-3}^{-1}A_{1,s+6k}^{-1},
\quad
M_{\ell+1}=M_{\ell}A_{2,s+6k-9}^{-1},
\\
M_{\ell+2}=M_{\ell+1}A_{2,s+6k-15}^{-1},
\quad
\ldots,
\quad
M_{k+\ell-1}
=M_{k+\ell-2}A_{2,s+3}^{-1}.
\end{gather*}
The dominant monomials in $\chi_q\Big(\mathcal{B}_{k, \ell}^{(s)}\Big)\chi_q\Big(\mathcal{B}_{k-1, \ell-1}^{(s+6)}\Big)$
are $M_0,\ldots,M_{k+\ell-2}$.

\item[$(2)$] Let $M=C_{k, \ell-1}^{(s)}C_{k-1, \ell}^{(s+6)}$, $k\geq 1$, $\ell \geq 1$.
Then dominant monomials
in $\chi_q\Big(\mathcal{C}_{k, \ell-1}^{(s)}\Big) \chi_q\Big(\mathcal{C}_{k-1, \ell}^{(s+6)}\Big)$ are
\begin{gather*}
M_0=M,
\quad
M_1=MA_{2,s+6k+6\ell-5}^{-1},
\quad
M_2=M_1A_{2,s+6k+6\ell-11}^{-1},
\quad
\ldots,
\\
M_{\ell-1}=M_{\ell-2}A_{2,s+6k+7}^{-1},
\quad
M_{\ell}=M_{\ell-1}A_{2,s+6k-3}^{-1}A_{1,s+6k}^{-1}A_{1,s+6k-2}
^{-1}A_{2,s+6k+1}^{-1},
\\
M_{\ell+1}=M_{\ell}A_{2,s+6k-9}^{-1},
\quad
M_{\ell+2}=M_{\ell+1}A_{2,s+6k-15}^{-1},
\quad
\ldots,
\quad
M_{k+\ell-1}
=M_{k+\ell-2}A_{2,s+3}^{-1}.
\end{gather*}
The dominant monomials in $\chi_q\Big(\mathcal{C}_{k, \ell}^{(s)}\Big)\chi_q\Big(\mathcal{C}_{k-1, \ell-1}^{(s+6)}\Big)$
are $M_0,\ldots,M_{k+\ell-2}$.

\item[$(3)$] Let $M=D_{0, \ell-1}^{(s)}B_{\ell, 0}^{(s+8)}$, $\ell \geq 1$.
Then dominant monomials in $\chi_q\Big(\mathcal{D}_{0, \ell-1}^{(s)}\Big)\chi_q\Big(\mathcal{B}_{\ell, 0}^{(s+8)}\Big)$
are
\begin{gather*}
M_0=M,
\quad
M_1=MA_{2,s+6\ell-1}^{-1},
\quad
M_2=M_1A_{2,s+6\ell-7}^{-1},
\quad
\ldots,
\\
M_{\ell-1}=M_{\ell-2}A_{2,s+11}^{-1},
\quad
M_{\ell}=M_{\ell-1}A_{1,s+6k+2}^{-1}A_{2,s+6k+5}^{-1}.
\end{gather*}
The dominant monomials in $\chi_q( \mathcal{D}_{0, \ell}^{(s)} ) \chi_q(\mathcal{B}_{\ell-1, 0}^{(s+8)} )$
are $M_0,\ldots,M_{\ell-1}$.

\item[$(4)$] Let $M{=}D_{k, \ell-1}^{(s)}D_{k-1, \ell}^{(s+6)}$, $k\geq 1$, $\ell \geq 1$.
Then dominant monomials
in $\chi_q\Big(\mathcal{D}_{k, \ell-1}^{(s)}\Big)\chi_q\Big(\mathcal{D}_{k-1, \ell}^{(s+6)}\Big)$
are
\begin{gather*}
M_0=M,
\quad
M_1=MA_{2,s+6k+6\ell-1}^{-1},
\\
M_2=M_1A_{2,s+6k+6\ell-7}^{-1},
\quad
\ldots,
\quad
M_{\ell-1}=M_{\ell-2}A_{2,s+6k+11}^{-1},
\\
M_{\ell}=M_{\ell-1}A_{1,s+6k+2}^{-1}A_{2,s+6k+5}^{-1},
\quad
M_{\ell+1}=M_{\ell}A_{2,s+6k-3}^{-1}A_{1,s+6k}^{-1},
\\
M_{\ell+2}=M_{\ell+1}A_{2,s+6k-9}^{-1},
\quad
\ldots,
\quad
M_{k+\ell}=M_{k+\ell-1}A_{2,s+3}^{-1}.
\end{gather*}
The dominant monomials
in $\chi_q\Big(\mathcal{D}_{k, \ell}^{(s)}\Big)\chi_q\Big(\mathcal{D}_{k-1, \ell-1}^{(s+6)}\Big)$
are $M_0,\ldots,M_{k+\ell-1}$.

\item[$(5)$] Let $M=E_{k, \ell-1}^{(s)}E_{k-1, \ell}^{(s+2)}$, $k\geq 1$, $\ell \geq 1$.

If $\ell=2r+1$, then dominant monomials
in $\chi_q\Big(\mathcal{E}_{k, \ell-1}^{(s)}\Big)\chi_q\Big(\mathcal{E}_{k-1, \ell}^{(s+2)}\Big)$ are
\begin{gather*}
M_0=M,
\quad
M_1=MA_{2,s+2k+3\ell-3}^{-1},
\\
M_2=M_1A_{2,s+2k+3\ell-9}^{-1},
\quad
\ldots,
\quad
M_{r}=M_{r-1}A_{2,s+2k+6}^{-1},
\\
M_{r+1}=M_{r}A_{1,s+2k-1}^{-1}A_{1,s+2k-3}^{-1}A_{2,s+2k}^{-1},
\quad
M_{r+2}=M_{r+1}A_{1,s+2k-5}^{-1},
\\
M_{r+3}=M_{r+2}A_{1,s+2k-7}^{-1},
\quad
\ldots,\quad M_{k+r-1}=M_{k+r-2}A_{1,s+1}
^{-1}.
\end{gather*}
The dominant monomials in $\chi_q\Big(\mathcal{E}_{k, \ell}^{(s)}\Big)\chi_q\Big(\mathcal{E}_{k-1, \ell-1}^{(s+2)}\Big)$
are $M_0,\ldots,M_{k+r-2}$.

If $\ell=2r$, then dominant monomials
in $\chi_q\Big(\mathcal{E}_{k, \ell-1}^{(s)}\Big)\chi_q\Big(\mathcal{E}_{k-1,\ell}^{(s+2)}\Big)$ are
\begin{gather*}
M_0=M,
\quad
M_1=MA_{2,s+2k+3\ell-4}^{-1},
\\
M_2=M_1A_{2,s+2k+3\ell-10}^{-1},
\quad
\ldots,
\quad
M_{r-1}=M_{r-2}A_{2,s+2k+8}^{-1},
\\
M_{r}=M_{r-1}A_{1,s+2k-1}^{-1}A_{2,s+2k+2}^{-1},
\quad
M_{r+1}=M_{r}A_{1,s+2k-3}^{-1},
\\
M_{r+2}=M_{r+1}A_{1,s+2k-5}^{-1},
\quad
\ldots,\quad M_{k+r-1}=M_{k+r-2}A_{1,s+1}
^{-1}.
\end{gather*}
The dominant monomials in $\chi_q\Big(\mathcal{E}_{k, \ell}^{(s)}\Big)\chi_q\Big(\mathcal{E}_{k-1, \ell-1}^{(s+2)}\Big)$
are $M_0,\ldots,M_{k+r-2}$.

\item[$(6)$] Let $M=F_{k, \ell-1}^{(s)}F_{k-1, \ell}^{(s+2)}$, $k\geq 1$, $\ell \geq 1$.
Then dominant monomials
in $\chi_q\Big(\mathcal{F}_{k, \ell-1}^{(s)}\Big)\chi_q\Big(\mathcal{F}_{k-1, \ell}^{(s+2)}\Big)$ are
\begin{gather*}
M_0=M,
\quad
M_1=MA_{1,s+2k+2\ell+3}^{-1},
\\
M_2=M_1A_{1,s+2k+2\ell+1}^{-1},
\quad
\ldots,
\quad
M_{\ell-1}=M_{\ell-2}A_{1,s+2k+7}^{-1},
\\
M_{\ell}=M_{\ell-1}A_{1,s+2k-1}^{-1}A_{2,s+2k+2}^{-1}A_{1,s+2k+5}^{-1},
\quad
M_{\ell+1}=M_{\ell}A_{1,s+2k-3}^{-1},
\\
M_{\ell+2}=M_{\ell+1}A_{1,s+2k-5}^{-1},
\quad
\ldots,
\quad
M_{k+\ell-1}
=M_{k+\ell-2}A_{1,s+1}^{-1}.
\end{gather*}
The dominant monomials in $\chi_q\Big(\mathcal{F}_{k, \ell}^{(s)}\Big)\chi_q\Big(\mathcal{F}_{k-1, \ell-1}^{(s+2)}\Big)$
are $M_0,\ldots,M_{k+\ell-2}$.
\end{enumerate}
In each case, for each $i$, the multiplicity of $M_i$ in the corresponding product of $q$-characters is~$1$.
\end{Lemma}

\begin{proof}
We prove the case of $\chi_q\Big(\mathcal{C}_{k, \ell-1}^{(s)}\Big)\chi_q\Big(\mathcal{C}_{k-1, \ell}^{(s+6)}\Big)$.
The other cases are similar.
Let $m'_1=C_{k, \ell-1}^{(s)}$, $m'_2=C_{k-1, \ell}^{(s+6)}$.
Without loss of generality, we assume that $s=6$.
Then
\begin{gather*}
m'_1=(2_{6}2_{12}\cdots2_{6k})(2_{6k+10}2_{6k+16}\cdots2_{6k+6\ell-2}),
\\
m'_2=(2_{12}\cdots2_{6k})(2_{6k+10}2_{6k+16}\cdots2_{6k+6\ell-2}2_{6k+6\ell+4}).
\end{gather*}

Let $m=m_1m_2$ be a~dominant monomial, where $m_i \in \chi_q(m'_i)$, $i=1, 2$.
Denote by $m_3=2_{6k+10}2_{6k+16}\cdots 2_{6k+6\ell +4}$.
If $m_2 \in \chi_q(2_{12}\cdots 2_{6k})(\chi_q(m_3)-m_3)$, then $m=m_1m_2$ is right negative and hence $m$
is not dominant.
Therefore $m_2 \in \chi_q(2_{12}\cdots 2_{6k})m_3$.

Suppose that $m_2 \in \mathscr{M}(L(m'_2)) \cap \mathscr{M}(\chi_q(2_{12}\cdots
2_{6k-6})(\chi_q(2_{6k})-2_{6k})m_3)$.
By the FM algorithm for $L(m'_2)$ and Lemma~\ref{fundamental modules}, $m_2$ must have a~factor
$2_{6k+6}^{-1}$ or $1_{6k+7}^{-1}$ or $2_{6k+8}^{-1}$.
By Lemma~\ref{fundamental modules}, $m_1$ does not have the factors $2_{6k+6}$ and $2_{6k+8}$.
Therefore $m_2$ cannot have factors $2_{6k+6}^{-1}$ and $2_{6k+8}^{-1}$ since $m=m_1m_2$ is dominant.
Hence $1_{6k+7}^{-1}$ is a~factor of $m_2$.
Since $m=m_1m_2$ is dominant, we need to cancel $1_{6k+7}^{-1}$ using a~factor in $m_1$.
By Lemma~\ref{fundamental modules}, the only possible way to cancel $1_{6k+7}^{-1}$ by $m_1$ is to use the
factor $1_{6k+5}1_{6k+7}1_{6k+9}2_{6k+8}^{-1}2_{6k+10}^{-1}$ or
$1_{6k+5}1_{6k+7}1_{6k+11}^{-1}2_{6k+8}^{-1}$ of $m_1$ coming from $\chi_q(2_{6k})$.
Since $2_{6k+8}^{-1}$ cannot be canceled by any monomials in $\chi_q(2_{6}2_{12}\cdots 2_{6k-6})$, we have
the factor $2_{6k+8}^{-1}$ in $m=m_1m_2$ and hence $m$ is not dominant.
Therefore $m_2 \in \mathscr{M}(L(2_{12}\cdots 2_{6k-6}))2_{6k}m_3$.
By the FM algorithm, $m_2=m'_2$.

If
\begin{gather*}
m_1\in\chi_q(2_{6}\cdots2_{6k}2_{6k+10}2_{6k+16}\cdots2_{6k+6\ell-8})(\chi_q(2_{6k+6\ell-2})
\\
\hphantom{m_1\in}{}
-2_{6k+6\ell-2}-2_{6k+6\ell+4}^{-1}1_{6k+6\ell-1}1_{6k+6\ell+1}1_{6k+6\ell+3}),
\end{gather*}
then $m=m_1m_2$ is right-negative and hence not dominant.
Therefore $m_1$ is in one of the following polynomials
\begin{gather}
\chi_q(2_{6}\cdots2_{6k}2_{6k+10}2_{6k+16}\cdots2_{6k+6\ell-8})2_{6k+6\ell-2},
\label{first set}
\\
\chi_q(2_{6}\cdots2_{6k}2_{6k+10}2_{6k+16}\cdots2_{6k+6\ell-8})2_{6k+6\ell+4}^{-1}1_{6k+6\ell-1}
1_{6k+6\ell+1}1_{6k+6\ell+3}.
\label{second set}
\end{gather}

If $m_1$ is in~\eqref{first set}, then $m_1=m'_1$.
The dominant monomial we obtain is $M_0=m'_1m'_2$.
If $m_1$ is the highest monomial in~\eqref{second set}, then we obtain the dominant monomial $M_1=m_1m'_2$.
Suppose that $m_1$ is in
\begin{gather*}
\begin{split}
& \mathscr{M}(L(m'_1))
\cap\mathscr{M}(\chi_q(2_{6}\cdots2_{6k}2_{6k+10}2_{6k+16}
\cdots2_{6k+6\ell-14})(\chi_q(2_{6k+6\ell-8})-2_{6k+6\ell-8})
\\
& \hphantom{\mathscr{M}(L(m'_1))
\cap}{}
\times
2_{6k+6\ell+4}^{-1}1_{6k+6\ell-1}1_{6k+6\ell+1}1_{6k+6\ell+3}).
\end{split}
\end{gather*}
By the FM algorithm for $L(m'_1)$,
\begin{gather*}
m_1\in\chi_q(2_{6}\cdots2_{6k}2_{6k+10}2_{6k+16}\cdots2_{6k+6\ell-14})
\\
\hphantom{m_1\in}{}
\times
\big(2_{6k+6\ell-2}^{-1}1_{6k+6\ell-7}1_{6k+6\ell-5}1_{6k+6\ell-3}\big)
\big(2_{6k+6\ell+4}^{-1}1_{6k+6\ell-1}1_{6k+6\ell+1}1_{6k+6\ell+3}\big).
\end{gather*}
We obtain the dominant monomial $M_2=m_1m'_2$.
Continue this procedure, we obtain dominant monomials $M_3,\ldots,M_{\ell-1}$ and the remaining dominant
monomials are of the form $m_1m'_2$, where $m_1$ is a~non-highest monomial in
\begin{gather*}
\mathscr{M}(L(m'_1))\cap\mathscr{M}(L(2_{6}\cdots2_{6k}))2_{6k+16}^{-1}2_{6k+22}^{-1}\cdots2_{6k+6\ell+4}
^{-1}1_{6k+11}1_{6k+13}\cdots1_{6k+6\ell+3}.
\end{gather*}
Suppose that $m_1$ is a~non-highest monomial in the above set.
Since the non-highest monomials in $\chi_q(2_{6}\cdots 2_{6k})$ are right-negative, we need cancellations of
factors with negative powers of some monomial in $\chi_q(2_{6}\cdots 2_{6k})$ with
$2_{6k+10}1_{6k+11}1_{6k+13} \cdots 1_{6k+6\ell+3}$.
The only cancellation can happen is to cancel $2_{6k+10}$ or $1_{6k+11}$.
Since $1_{6k+9}^{2}$ does not appear in $\chi_q(2_{6}\cdots 2_{6k})$, $1_{6k+11}$ cannot be canceled.
Therefore we need a~cancellation with $2_{6k+10}$.
The only monomials in $\chi_q(2_{6}\cdots 2_{6k})$ which can cancel $2_{6k+10}$ is in one of the following
polynomials
\begin{gather*}
\chi_q(2_{6}\cdots2_{6k-6})1_{6k+1}1_{6k+9}2_{6k+10}^{-1},
\\
\chi_q(2_{6}\cdots2_{6k-6})1_{6k+3}^{-1}1_{6k+9}2_{6k+2}2_{6k+10}^{-1},
\\
\chi_q(2_{6}\cdots2_{6k-6})1_{6k+5}1_{6k+7}1_{6k+9}2_{6k+8}^{-1}2_{6k+10}^{-1}.
\end{gather*}
Therefore $m_1$ is in one of the following sets
\begin{gather}
\mathscr{M}(L(m'_1))\cap\mathscr{M}(L(2_{6}\cdots2_{6k-6}))1_{6k+1}1_{6k+9}2_{6k+10}^{-1}
\cdots2_{6k+6\ell+4}^{-1}1_{6k+11}\cdots1_{6k+6\ell+3},
\label{one}
\\
\mathscr{M}(L(m'_1))\cap\mathscr{M}(L(2_{6}\cdots2_{6k-6}))
\nonumber
\\
\hphantom{\mathscr{M}(L(m'_1))\cap}{}
\times
1_{6k+3}^{-1}1_{6k+9}2_{6k+2}2_{6k+10}^{-1}\cdots2_{6k+6\ell+4}^{-1}1_{6k+11}\cdots1_{6k+6\ell+3},
\label{two}
\\
\mathscr{M}(L(m'_1))\cap\mathscr{M}(L(2_{6}\cdots2_{6k-6}))
\nonumber
\\
\hphantom{\mathscr{M}(L(m'_1))\cap}{}
\times
1_{6k+5}1_{6k+7}1_{6k+9}2_{6k+8}^{-1}2_{6k+10}^{-1}\cdots2_{6k+6\ell+4}^{-1}1_{6k+11}\cdots1_{6k+6\ell+3}.
\label{three}
\end{gather}

If $m_1$ is in~\eqref{two}, then we need to cancel $1_{6k+3}^{-1}$.
We have
\begin{gather*}
\mathscr{M}(L(2_{6}\cdots2_{6k-6}))\subset\mathscr{M}(\chi_q(2_{6}\cdots2_{6k-12})\chi_q(2_{6k-6})).
\end{gather*}
By Lemma~\ref{fundamental modules}, only the monomials
\begin{gather*}
1_{6k-5}1_{6k+3}2_{6k+4}^{-1},
\qquad
1_{6k-3}^{-1}1_{6k+3}2_{6k-4}2_{6k+4}^{-1},
\qquad
1_{6k-1}1_{6k+1}1_{6k+3}
2_{6k+2}^{-1}2_{6k+4}^{-1}
\end{gather*}
in $\chi_q(2_{6k-6})$ can cancel $1_{6k+3}^{-1}$.
But these monomials have the factor $2_{6k+4}^{-1}$ which cannot be canceled by any monomials in
$\chi_q(2_{6}\cdots 2_{6k-12})$ or by $m'_2$.
Hence $m_1$ is not in~\eqref{two}.

If $m_1$ is in~\eqref{three}, then we need to cancel $2_{6k+8}^{-1}$.
But $2_{6k+8}^{-1}$ cannot be canceled by any monomials in $\chi_q(2_{6}\cdots 2_{6k-6})$ or by $m'_2$.
Therefore $m_1$ is not in~\eqref{three}.
Hence $m_1$ is in~\eqref{one}.

If $m_1$ is the highest monomial in~\eqref{one} with respect to $\leq$ def\/ined in~\eqref{partial order of
monomials}, then $m_1m'_2=M_{\ell}$.
Suppose that $m_1$ a~non-highest monomial in~\eqref{one}.
By the FM algorithm, $m_1$~must be in
\begin{gather*}
\chi_q(2_{6}\cdots 2_{6k-12})2_{6k}^{-1}1_{6k-5}1_{6k-3}1_{6k-1}1_{6k+1}1_{6k+9}2_{6k+10}^{-1}
\cdots 2_{6k+6\ell+4}^{-1}1_{6k+11}\cdots 1_{6k+6\ell+3}.
\end{gather*}
If $m_1$ is the highest monomial in the above set, then $m_1m'_2=M_{\ell+1}$.
Continue this procedure, we can show that the only remaining dominant monomials are $M_{\ell+2},\ldots,M_{k+\ell-1}$.

It is clear that the multiplicity of $M_i$, $i=1,\ldots,k+\ell-1$, in $\chi_q(m_1) \chi_q(m_2)$ is~$1$.
\end{proof}

\subsection{Products of sources are special}
\begin{Lemma}
\label{dominant monomials in sources}
Let $[\mathcal{S}]$ be the last summand in one of the relations~\eqref{b equation}--\eqref{f equation}.
Then $\mathcal{S}$ is special.
\end{Lemma}
\begin{proof}
We give a~proof for $\mathcal{S}$ in the last line of~\eqref{e equation} and in the last line of~\eqref{f
equation}.
The other cases are similar.

Let $S_1=\chi_q( \mathcal{C}_{r+1, p}^{(s+1)} ) \chi_q( \mathcal{B}_{r+1, 3p-1}^{(s+3)} ) \chi_q(
\mathcal{B}_{r+p+1, 0}^{(s+5)})$.
Let
\begin{gather*}
n_1=2_{s+1}2_{s+7}\cdots2_{s+6r-5}2_{s+6r+1},
\qquad
n'_1=2_{s+6r+11}2_{s+6r+17}\cdots2_{s+6r+6p+5},
\\
n_2=2_{s+3}2_{s+9}\cdots2_{s+6r-3}2_{s+6r+3},
\qquad
n'_2=1_{s+6r+10}1_{s+6r+12}\cdots1_{s+6r+6p+6},
\\
n_3=2_{s+5}2_{s+11}\cdots2_{s+6r+6p+5}.
\end{gather*}
Then $C_{r+1, p}^{(s+1)}=n_1n'_1$, $B_{r+1, 3p-1}^{(s+3)}=n_2n'_2$, $B_{r+p+1, 0}^{(s+5)}=n_3$.
Let $m'=m_1m_2m_3$ be a~dominant monomial, where
\begin{gather*}
m_1\in\mathscr{M}\Big(\mathcal{C}_{r+1,p}^{(s+1)}\Big),
\qquad
m_2\in\mathscr{M}\Big(\mathcal{B}_{r+1,3p-1}^{(s+3)}
\Big),
\qquad
m_3\in\mathscr{M}\Big(\mathcal{B}_{r+p+1,0}^{(s+5)}\Big).
\end{gather*}

If $m_3 \neq B_{r+p+1, 0}^{(s+5)}$ or $m_1 \in \chi_q(n_1)(\chi_q(n'_1) - n'_1 )$ or $m_2 \in
\chi_q(n_2)(\chi_q(n'_2) - n'_2 )$, then $m'$ is right-negative which contradicts the fact that $m'$ is
dominant.
Therefore $m_3= B_{r+p+1, 0}^{(s+5)}$, $m_1 \in \chi_q(n_1)n'_1$, and $m_2 \in \chi_q(n_2)n'_2$.

If $m_2$ is in
\begin{gather}
\mathscr{M}(L(n_2n'_2))\cap\mathscr{M}(\chi_q(2_{s+3}2_{s+9}\cdots2_{s+6r-3})(\chi_q(2_{s+6r+3})-2_{s+6r+3}
)n'_2),
\label{a set}
\end{gather}
then
\begin{gather*}
m_2\in\chi_q(2_{s+3}2_{s+9}\cdots2_{s+6r-3})2_{s+6r+9}^{-1}1_{s+6k+4}1_{s+6k+6}1_{s+6k+8}n'_2.
\end{gather*}
By Lemma~\ref{fundamental modules}, the factor $2_{s+6r+9}^{-1}$ cannot be canceled by any monomial in
either $\chi_q(n_1)$ or $\chi_q(2_{s+3}2_{s+9}\cdots 2_{s+6r-3})$.
It is clear that $2_{s+6r+9}^{-1}$ cannot be canceled by $n'_1$, $n'_2$, $n_3$.
Therefore $2_{s+6r+9}^{-1}$ cannot be canceled.
Hence $m_2$ is not in~\eqref{a set}.
Thus $m_2$ must be in
\begin{gather*}
\mathscr{M}(L(n_2n'_2))\cap\mathscr{M}(L(2_{s+3}2_{s+9}\cdots2_{s+6r-3}))2_{s+6r+3}n'_2.
\end{gather*}
Therefore $m_2=B_{r+1, 3p-1}^{(s+3)}$.

Suppose that $m_1 \neq C_{r+1, p}^{(s+1)}$.
Then $m_1=m'_1n'_1$, where $m'_1$ is a~non-highest monomial in~$\chi_q(n_1)$.
Since the non-highest monomials in $\chi_q(n_1)$ are right-negative, we need a~cancellation with
$n'_1n'_2m_3$.
The only cancellation can happen is to cancel $2_{s+6r+11}$ in $n'_1$, or cancel one of $2_{s+6r+3}$,
$1_{s+6r+10}$ in $n_2n'_2$, or cancel one of $2_{s+6r+5}$, $2_{s+6r+11}$ in $m_3$.
By the FM algorithm, $2_{s+6r+11}$ cannot be canceled.
By Lemma~\ref{fundamental modules}, $1_{s+6r+10}$, $2_{s+6r+3}$ and $2_{s+6r+5}$ cannot be canceled.
This is a~contradiction.
Therefore $m_1= C_{r+1, p}^{(s+1)}$.

Therefore the only dominant monomial in $S_1$ is $C_{r+1, p}^{(s+1)} B_{r+1, 3p-1}^{(s+3)} B_{r+p+1,
0}^{(s+5)}$.

Let $S_2=\chi_q\Big(\mathcal{D}_{r+1, \left\lfloor \frac{\ell}{3} \right\rfloor}^{(s+1)}\Big)
\chi_q\Big(\mathcal{C}_{r+1,\left\lfloor \frac{\ell+1}{3} \right\rfloor}^{(s+3)}\Big)
\chi_q\Big(\mathcal{B}_{r, 0}^{(s+5)}\Big)
\chi_q\Big(\mathcal{B}_{\left\lfloor \frac{\ell-1}{3} \right\rfloor, 0}^{(s+6r+17)}\Big)$,
$r\geq 0$, and $\ell = 3p$, $p\geq 1$.
The cases of $\ell=3p+1$, $p\geq 0$ and $\ell=3p+2$, $p\geq 0$ are similar.
Let
\begin{gather*}
n_1=2_{s+1}2_{s+7}\cdots2_{s+6r+1}1_{s+6r+8},
\qquad
n'_1=2_{s+6r+15}2_{s+6r+21}\cdots2_{s+6r+6p+9},
\\
n_2=2_{s+3}2_{s+9}\cdots2_{s+6r-3}2_{s+6r+3},
\qquad
n'_2=2_{s+6r+13}2_{s+6r+20}\cdots2_{s+6r+6p+7},
\\
n_3=2_{s+5}2_{s+11}\cdots2_{s+6r-1},
\qquad
n_4=2_{s+6r+17}2_{s+6r+23}\cdots2_{s+6r+6p+5}.
\end{gather*}
Then $D_{r+1, \left\lfloor \frac{\ell}{3} \right\rfloor}^{(s+1)}
= n_1n'_1$,
$C_{r+1, \left\lfloor \frac{\ell+1}{3}\right\rfloor}^{(s+3)} = n_2n'_2$,
$B_{r, 0}^{(s+5)} = n_3$,
$B_{\left\lfloor \frac{\ell-1}{3} \right\rfloor,0}^{(s+6r+17)} = n_4$.

Let $m'=m_1m_2m_3m_4$ be a~dominant monomial, where
\begin{gather*}
m_1 \in \mathscr{M}\Big(\mathcal{D}_{r+1,\left\lfloor\frac{\ell}{3}\right\rfloor}^{(s+1)}\Big),
\qquad
m_2 \in \mathscr{M}\Big(\mathcal{C}_{r+1,\left\lfloor\frac{\ell+1}{3}\right\rfloor}^{(s+3)}\Big),
\\
m_3 \in \mathscr{M}\Big(\mathcal{B}_{r,0}^{(s+5)}\Big),
\qquad
m_4 \in \mathscr{M}\Big(\mathcal{B}_{\left\lfloor\frac{\ell-1}{3}\right\rfloor,0}^{(s+6r+17)}\Big).
\end{gather*}

If $m_4 \neq n_4$ or $m_1 \in \chi_q(n_1)(\chi_q(n'_1) - n'_1 )$ or $m_2 \in \chi_q(n_2)(\chi_q(n'_2) -
n'_2)$, then $m'$ is right-negative which contradicts the fact that $m'$ is dominant.
Therefore $m_4= n_4$, $m_1 \in \chi_q(n_1)n'_1$, and $m_2 \in \chi_q(n_2)n'_2$.

If
\begin{gather}
m_1\in\mathscr{M}(L(n_1n'_1))\cap\mathscr{M}(\chi_q(2_{s+1}2_{s+7}\cdots2_{s+6r+1})(\chi_q(1_{s+6r+8}
)-1_{s+6r+8})n'_1),
\end{gather}
then by the FM algorithm for $L(n_1n'_1)$,
\begin{gather*}
m_1\in\chi_q(2_{s+1}2_{s+7}\cdots2_{s+6r+1})1_{s+6r+10}^{-1}2_{s+6r+9}n'_1.
\end{gather*}
It is clear that $1_{s+6r+10}^{-1}$ is not canceled by $n'_1$, $n'_2$, $n_4$, and any monomial in
$\chi_q(n_3)$.
By the FM algorithm for $\chi_q(n_2n'_2)$, $1_{s+6r+10}^{-1}$ cannot be canceled by any monomial in
$\chi_q(n_2n'_2)$.
Therefore, by Lemma~\ref{fundamental modules}, $1_{s+6r+10}^{-1}$ can only be canceled by one of the factors
\begin{gather*}
1_{s+6r+2}1_{s+6r+10}2_{s+6r+11}^{-1},
\\
1_{s+6r+4}^{-1}1_{s+6r+10}2_{s+6r+3}2_{s+6r+11}^{-1},
\qquad
1_{s+6r+6}1_{s+6r+8}1_{s+6r+10}2_{s+6r+9}^{-1}
2_{s+6r+11}^{-1}
\end{gather*}
coming from $\chi_q(2_{s+6r+1})$, where $2_{s+6r+1}$ is in $n_1$.
But then $2_{s+6r+11}^{-1}$ cannot be canceled.
This contradicts the fact that $m'$ is dominant.
Hence $m_1$ is not in~\eqref{a set}.
Thus $m_1$ must be in
\begin{gather*}
\mathscr{M}(L(n_1n'_1))\cap\mathscr{M}(L(n_1n'_1))\cap\mathscr{M}(L(2_{s+1}2_{s+7}\cdots2_{s+6r+1}
))1_{s+6r+8}n'_1.
\end{gather*}
If $m_1$ is in
\begin{gather*}
\mathscr{M}(L(n_1n'_1))\cap\mathscr{M}(L(n_1n'_1))\cap\mathscr{M}(\chi_q(2_{s+1}
\\
\hphantom{\mathscr{M}(L(n_1n'_1))\cap}{}
\times2_{s+7}\cdots 2_{s+6r-5})(\chi_q(2_{s+6r+1})-2_{s+6r+1})1_{s+6r+8}n'_1).
\end{gather*}
Then
\begin{gather*}
m_1\in\chi_q(2_{s+1}2_{s+7}\cdots2_{s+6r-5})2_{s+6r+7}^{-1}1_{s+6r+2}1_{s+6r+4}1_{s+6r+6}1_{s+6r+8}n'_1.
\end{gather*}
The only possible way to cancel $2_{s+6r+7}^{-1}$ is to use one of the terms
\begin{gather*}
1_{s+6r+4}1_{s+6r+8}^{-1}1_{s+6r+10}^{-1}2_{s+6r+7},
\\
1_{s+6r+6}^{-1}1_{s+6r+8}^{-1}1_{s+6r+10}^{-1}
2_{s+6r+5}2_{s+6r+7},
\qquad
2_{s+6r+7}2_{s+6r+11}^{-1}
\end{gather*}
in $\chi_q(2_{s+6r+3})$, where $2_{s+6r+3}$ is in $n_2$.
But then we have to cancel $1_{s+6r+10}^{-1}$ or $2_{s+6r+11}^{-1}$.
But $1_{s+6r+10}^{-1}$ and $2_{s+6r+11}^{-1}$ cannot be canceled.
This is a~contradiction.
Therefore $m_1$ must be in
\begin{gather*}
\mathscr{M}(L(n_1n'_1))\cap\mathscr{M}(L(2_{s+1}2_{s+7}\cdots2_{s+6r-5}))2_{s+6r+1}1_{s+6r+8}n'_1.
\end{gather*}
Hence $m_1=n_1n'_1$.

By the FM algorithm, when we compute the $q$-character for $\chi_q(n_2n'_2)$, we can only choose one of the
following terms
\begin{gather*}
2_{s+6r+3},
\qquad
1_{s+6r+4}1_{s+6r+6}1_{s+6r+8}2_{s+6r+9}^{-1},
\qquad
1_{s+6r+4}1_{s+6r+6}1_{s+6r+10}^{-1},
\\
1_{s+6r+4}1_{s+6r+8}^{-1}1_{s+6r+10}^{-1}2_{s+6r+7},
\\
1_{s+6r+6}^{-1}1_{s+6r+8}^{-1}1_{s+6r+10}^{-1}2_{s+6r+5}2_{s+6r+7},
\qquad
2_{s+6r+11}^{-1}2_{s+6r+7}
\end{gather*}
in $\chi_q(2_{s+6r+3})$.
Since $2_{s+6r+9}^{-1}$, $1_{s+6r+10}^{-1}$, and $2_{s+6r+11}^{-1}$ cannot be canceled, we can only choose
$2_{s+6r+3}$.
Therefore $m_2$ is in
\begin{gather*}
\mathscr{M}(L(n_2n'_2))\cap\mathscr{M}(L(2_{s+3}2_{s+9}\cdots2_{s+6r-3}))2_{s+6r+3}n'_2.
\end{gather*}
Therefore $m_2=n_2n'_2$.

If $m_3$ is in
\begin{gather*}
\mathscr{M}(L(n_3))\cap\mathscr{M}(\chi_q(2_{s+5}2_{s+11}\cdots2_{s+6r-7})(\chi_q(2_{s+6r-1})-2_{s+6r-1})),
\end{gather*}
then, by Lemma~\ref{fundamental modules}, $m=m_1m_2m_3m_4$ is non-dominant since $m_1=n_1n'_1, m_2=n_2n'_2,
m_4=n_4$.
This contradicts the fact that $m$ is dominant.
Therefore $m_3$ is in
\begin{gather*}
\mathscr{M}(L(n_3))\cap\mathscr{M}(L(2_{s+5}2_{s+11}\cdots2_{s+6r-7}))2_{s+6r-1}.
\end{gather*}
Hence $m_3=n_3$.

Therefore the only dominant monomial in $S_2$ is $D_{r+1, \left\lfloor \frac{\ell}{3} \right\rfloor}^{(s+1)}
C_{r+1, \left\lfloor \frac{\ell+1}{3} \right\rfloor}^{(s+3)} B_{r, 0}^{(s+5)} B_{\left\lfloor \frac{\ell-1}{3}
\right\rfloor, 0}^{(s+6r+17)}$.
\end{proof}

\subsection{Proof of Theorem~\ref{extended T-system}} By Lemmas~\ref{dominant monomials} and~\ref{dominant
monomials in sources}, the dominant monomials in the q-characters of the left hand side and of the right
hand side of every relation in Theorem~\ref{extended T-system} are the same.
The theorem follows.

\section{Proof of Theorem~\ref{irreducible}}
\label{proof irreducible}

By Lemma~\ref{dominant monomials in sources}, $\mathcal{S}$ is special and hence irreducible.
Therefore we only have to show that $\mathcal{T}\otimes \mathcal{B}$ is irreducible.
It suf\/f\/ices to prove that for each non-highest dominant monomial $M$ in $\mathcal{T} \otimes
\mathcal{B}$, we have $\mathscr{M}(L(M)) \not\subset \mathscr{M}(\mathcal{T} \otimes \mathcal{B})$.
The idea is similar as in~\cite{Her06,MY11b, Nak04}.
Recall that the dominant monomials in $\mathcal{T}\otimes \mathcal{B}$ are described by Lemma~\ref{dominant
monomials}.
\begin{Lemma}
We consider the same cases as in Lemma~{\rm \ref{dominant monomials}}.
In each case $M_i$ are the dominant monomials described by that lemma.
\begin{enumerate}\itemsep=0pt
\item[$(1)$] For $k\geq 1$, $\ell \geq 1$, let
\begin{gather*}
n_1=M_1A_{1,s+6k+2\ell-2}^{-1},
\quad
n_2=M_2A_{1,s+6k+2\ell-4}^{-1},
\quad
\ldots,
\\
n_{\ell-1}=M_{\ell-1}A_{1,s+6k+2}^{-1},
\quad
n_{\ell}=M_{\ell}A_{2,s+6k-3}^{-1}A_{1,s+6k}^{-1},
\\
n_{\ell+1}=M_{\ell+1}A_{2,s+6k-9}^{-1},
\quad
\ldots,
\quad
n_{k+\ell-2}=M_{k+\ell-2}A_{2,s+9}^{-1}.
\end{gather*}
Then for $i=1,\ldots,k+\ell-2$, $n_i\in \chi_q(M_i)$ and $n_i \not\in \chi_q\Big(\mathcal{B}_{k, \ell}^{(s)}\Big)
\chi_q\Big(\mathcal{B}_{k-1, \ell-1}^{(s+6)}\Big)$.

\item[$(2)$] For $k\geq 1$, $\ell \geq 1$, let
\begin{gather*}
n_1=M_1A_{2,s+6k+6\ell-5}^{-1},
\quad
n_2=M_2A_{2,s+6k+6\ell-11}^{-1},
\quad
\ldots,
\\
n_{\ell-1}=M_{\ell-1}A_{2,s+6k+7}^{-1},
\quad
n_{\ell}=M_{\ell}A_{2,s+6k-3}^{-1}A_{1,s+6k}^{-1}A_{1,s+6k-2}^{-1}
A_{2,s+6k+1}^{-1},
\\
n_{\ell+1}=M_{\ell+1}A_{2,s+6k-9}^{-1},
\quad
\ldots,
\quad
n_{k+\ell-2}=M_{k+\ell-2}A_{2,s+9}^{-1}.
\end{gather*}
Then for $i=1,\ldots,k+\ell-2$, $n_i\in \chi_q(M_i)$ and $n_i \not\in \chi_q\Big(\mathcal{C}_{k, \ell}^{(s)}\Big)
\chi_q\Big(\mathcal{C}_{k-1, \ell-1}^{(s+6)}\Big)$.

\item[$(3)$] For $\ell \geq 1$, let
\begin{gather*}
n_1=M_1A_{2,s+6\ell-1}^{-1},
\quad
n_2=M_2A_{2,s+6\ell-7}^{-1},
\quad
\ldots,\quad n_{\ell-1}=M_{\ell-1}A_{2,s+11}^{-1}.
\end{gather*}
Then for $i=1,\ldots,\ell$, $n_i\in \chi_q(M_i)$ and $n_i \not\in \chi_q\Big(\mathcal{D}_{0, \ell}^{(s)}\Big)
\chi_q\Big(\mathcal{B}_{\ell-1, 0}^{(s+8)}\Big)$.

\item[$(4)$] For $k\geq 1$, $\ell \geq 1$, let
\begin{gather*}
n_1=M_1A_{2,s+6k+6\ell-1}^{-1},
\quad
n_2=M_2A_{2,s+6k+6\ell-7}^{-1},
\quad
\ldots,
\quad
n_{\ell-1}=M_{\ell-1}A_{2,s+6k+11}^{-1},
\\
n_{\ell}=M_{\ell}A_{1,s+6k+2}^{-1}A_{2,s+6k+5}^{-1},
\quad
n_{\ell+1}=M_{\ell+1}A_{2,s+6k-3}^{-1}A_{1,s+6k}^{-1},
\\
n_{\ell+2}=M_{\ell+2}A_{2,s+6k-9}^{-1},
\quad
\ldots,
\quad
n_{k+\ell-1}=M_{k+\ell-1}A_{2,s+9}^{-1}.
\end{gather*}
Then for $i=1,\ldots,k+\ell-1$, $n_i\in \chi_q(M_i)$ and $n_i \not\in \chi_q\Big(\mathcal{D}_{k-1,\ell-1}^{(s+6)}\Big)
\chi_q\Big(\mathcal{D}_{k, \ell}^{(s)}\Big)$.

\item[$(5)$] For $k\geq 0, \ell=2r+1, r\geq 0$, let
\begin{gather*}
n_1=M_1A_{2,s+2k+3\ell-3}^{-1},
\quad
n_2=M_2A_{2,s+2k+3\ell-9}^{-1},
\quad
\ldots,
\\
n_{r}=M_{r}A_{2,s+2k+3}^{-1},
\quad
n_{r+1}=M_{r+1}A_{1,s+2k-1}^{-1}A_{1,s+2k-3}^{-1}A_{2,s+2k}^{-1},
\\
n_{r+2}=M_{r+2}A_{1,s+2k-5}^{-1},
\quad
\ldots,
\quad
n_{k+r-2}=M_{k+r-2}A_{1,s+3}^{-1}.
\end{gather*}
Then for $i=1,\ldots,r+k-2$, $n_i\in \chi_q(M_i)$ and $n_i \not\in \chi_q\Big(\mathcal{E}_{k, \ell}^{(s)}\Big)
\chi_q\Big(\mathcal{E}_{k-1, \ell-1}^{(s+2)}\Big)$.

For $k\geq 0, \ell=2r, r\geq 1$, let
\begin{gather*}
n_1=M_1A_{2,s+2k+3\ell-4}^{-1},
\quad
n_2=M_2A_{2,s+2k+3\ell-10}^{-1},
\quad
\ldots,
\\
n_{r-1}=M_{r-1}A_{2,s+2k+8}^{-1},
\quad
n_{r}=M_{r}A_{1,s+2k-1}^{-1}A_{2,s+2k+2}^{-1},
\\
n_{r+1}=M_{r+1}A_{1,s+2k-3}^{-1},
\quad
\ldots,
\quad
n_{k+r-2}=M_{k+r-2}A_{1,s+3}^{-1}.
\end{gather*}
Then for $i=1,\ldots,r+k-2$, $n_i\in \chi_q(M_i)$ and $n_i \not\in \chi_q\Big(\mathcal{E}_{k, \ell}^{(s)}\Big)
\chi_q\Big(\mathcal{E}_{k-1, \ell-1}^{(s+2)}\Big)$.

\item[$(6)$] For $k\geq 1$, $\ell \geq 1$, let
\begin{gather*}
n_1=M_1A_{1,s+2k+2\ell+3}^{-1},
\quad
n_2=M_2A_{1,s+2k+2\ell+1}^{-1},
\quad
\ldots,
\\
n_{\ell-1}=M_{\ell-1}A_{1,s+2k+7}^{-1},
\quad
n_{\ell}=M_{\ell}A_{1,s+2k-1}^{-1}A_{2,s+2k+2}^{-1}A_{1,s+2k+5}
^{-1},
\\
n_{\ell+1}=M_{\ell+1}A_{1,s+2k-3}^{-1},
\quad
\ldots,
\quad
n_{k+\ell-2}=M_{k+\ell-2}A_{1,s+3}^{-1}.
\end{gather*}
Then $i=1,\ldots,k+\ell-2$, $n_i\in \chi_q(M_i)$ and $n_i \not\in \chi_q\Big(\mathcal{F}_{k, \ell}^{(s)}\Big)
\chi_q\Big(\mathcal{F}_{k-1, \ell-1}^{(s+2)}\Big)$.
\end{enumerate}
\end{Lemma}

\begin{proof}
We give a~proof in the case of $\chi_q\Big(\mathcal{C}_{k, \ell}^{(s)}\Big)
\chi_q\Big(\mathcal{C}_{k-1,\ell-1}^{(s+6)}\Big)$.
The other cases are similar.
By def\/inition, we have
\begin{gather*}
C_{k,\ell}^{(s)}=(2_{s}2_{s+6}\cdots2_{s+6k-6})(2_{s+6k+4}2_{s+6k+10}\cdots2_{s+6k+6\ell-8}
2_{s+6k+6\ell-2}),
\\
C_{k-1,\ell-1}^{(s+6)}=(2_{s+6}2_{s+12}\cdots2_{s+6k-6})(2_{s+6k+4}2_{s+6k+10}\cdots2_{s+6k+6\ell-8}),
\\
M_1=C_{k,\ell}^{(s)}C_{k-1,\ell-1}^{(s+6)}A_{2,s+6k+6\ell-5}^{-1}
\\
\phantom{M_1}
=C_{k,\ell}^{(s)}C_{k-1,\ell-1}^{(s+6)}2_{s+6k+6\ell-8}^{-1}2_{s+6k+6\ell-2}^{-1}1_{s+6k+6\ell-7}
1_{s+6k+6\ell-5}1_{s+6k+6\ell-3}.
\end{gather*}
By $U_{q_2}(\hat{\mathfrak{sl}}_2)$ argument, it is clear that $n_1=M_1A_{2, s+6k+6\ell-5}^{-1}$ is in
$\chi_q(M_1)$.

If $n_1$ is in $\chi_q\Big(\mathcal{C}_{k, \ell}^{(s)}\Big)\chi_q\Big(\mathcal{C}_{k-1, \ell-1}^{(s+6)}\Big)$, then
$C_{k, \ell}^{(s)} A_{2, s+6k+6\ell-5}^{-1}$ is in $\chi_q\Big(\mathcal{C}_{k, \ell}^{(s)}\Big)$ which is
impossible by the FM algorithm for $\mathcal{C}_{k, \ell}^{(s)}$.
Similarly, $n_i\in \chi_q(M_i)$, $i=2,\ldots,\ell-1$, but $n_2,\ldots,n_{\ell-1}$
are not in $\chi_q\Big(\mathcal{C}_{k, \ell}^{(s)}\Big)\chi_q\Big(\mathcal{C}_{k-1, \ell-1}^{(s+6)}\Big)$.

By def\/inition,
\begin{gather*}
M_{\ell}=(2_{s}2_{s+6}\cdots2_{s+6k-6})(2_{s+6}2_{s+12}\cdots2_{s+6k-12})(1_{s+6k-5}1_{s+6k+3}1_{s+6k+5}
\cdots1_{s+6k+6\ell-3}).
\end{gather*}
Let $U=\big\{(1, aq^{s+6k}), (1, aq^{s+6k-3}), (2, aq^{s+6k-2}), (2, aq^{s+6k+1})\big\} \subset I \times
\mathbb{C}^{\times}$.
Let $\mathcal{M}$ be the f\/inite set consisting of the following monomials
\begin{gather*}
m_0=M_{\ell},
\qquad
m_1=m_0A_{2,s+6k-3}^{-1},
\qquad
m_2=m_1A_{1,s+6k}^{-1},
\\
m_{3}=m_{2}A_{1,s+6k-2}^{-1},
\qquad
m_4=m_{3}
A_{2,s+6k+1}^{-1}.
\end{gather*}
It is clear that $\mathcal{M}$ satisf\/ies the conditions in Theorem~\ref{truncated}.
Therefore
\begin{gather*}
\trunc_{m_+\mathcal{Q}_{U}^{-}}(\chi_q(M_{\ell}))=\sum_{m\in\mathcal{M}}m
\end{gather*}
and hence $n_{\ell}=M_{\ell}A_{2, s+6k-3}^{-1}A_{1, s+6k}^{-1}A_{1, s+6k-2}^{-1}A_{2, s+6k+1}^{-1}$ is in
$\chi_q(M_{\ell})$.

If $n_{\ell}$ is in $\chi_q\Big(\mathcal{C}_{k, \ell}^{(s)}\Big)\chi_q\Big(\mathcal{C}_{k-1, \ell-1}^{(s+6)}\Big)$, then
$C_{k, \ell}^{(s)} A_{2, s+6k-3}^{-1}A_{1, s+6k}^{-1}A_{1, s+6k-2}^{-1}A_{2, s+6k+1}^{-1}$ is
in $\chi_q\Big(\mathcal{C}_{k, \ell}^{(s)}\Big)$ which is impossible by the FM algorithm for $\mathcal{C}_{k,
\ell}^{(s)}$.
Similarly, we show that for $i=\ell+1,\ldots,k+\ell-2$, $n_i\in \chi_q(M_i)$ and $n_i \not\in
\chi_q\Big(\mathcal{C}_{k, \ell}^{(s)}\Big) \chi_q\Big( \mathcal{C}_{k-1, \ell-1}^{(s+6)}\Big)$.
\end{proof}

\section[The second part of the extended $T$-system]{The second part of the extended $\boldsymbol{T}$-system}
\label{The second part of the extended T-system}

Let $\tilde{B}_{k, \ell}^{(s)}$, $\tilde{C}_{k, \ell}^{(s)}$, $\tilde{D}_{k, \ell}^{(s)}$, $\tilde{E}_{k,
\ell}^{(s)}$, $\tilde{F}_{k, \ell}^{(s)}$ be the monomials obtained from $B_{k, \ell}^{(s)}$, $C_{k,
\ell}^{(s)}$, $D_{k, \ell}^{(s)}$, $E_{k, \ell}^{(s)}$, $F_{k, \ell}^{(s)}$ by replacing $i_a$ with
$i_{-a}$, $i=1, 2$.
Namely,
\begin{gather*}
\tilde{B}_{k,\ell}^{(s)}=\left(\prod_{i=0}^{\ell-1}1_{-s-6k-2i-1}\right)\left(\prod_{i=0}^{k-1}2_{-s-6i}
\right),
\qquad
\tilde{C}_{k,\ell}^{(s)}=\left(\prod_{i=0}^{\ell-1}2_{-s-6k-6i-4}\right)\left(\prod_{i=0}^{k-1}
2_{-s-6i}\right),
\\
\tilde{D}_{k,\ell}^{(s)}=\left(\prod_{i=0}^{\ell-1}2_{-s-6k-6i-8}\right)1_{-s-6k-1}\left(\prod_{i=0}^{k-1}
2_{-s-6i}\right),
\\
\tilde{F}_{k,\ell}^{(s)}=\left(\prod_{i=0}^{\ell-1}1_{-s-2k-2i-6}\right)\left(\prod_{i=0}
^{k-1}1_{-s-2i}\right),
\\
\tilde{E}_{k,\ell}^{(s)}=\left(\prod_{i=0}^{\left\lfloor\frac{\ell-2}{2}\right\rfloor}2_{-s-2k-6i-5}
\right)\left(\prod_{i=0}^{\left\lfloor\frac{\ell-1}{2}\right\rfloor}2_{-s-2k-6i-3}\right)\left(\prod_{i=0}^{k-1}
1_{-s-2i}\right).
\end{gather*}

Note that, in particular, for $k \in \mathbb{Z}_{\geq 0}$, $s\in \mathbb{Z}$, we have the following trivial
relations
\begin{gather}
\tilde{\mathcal{B}}_{k,0}^{(s)}=\tilde{\mathcal{C}}_{k,0}^{(s)}=\tilde{\mathcal{C}}_{0,k}^{(s-4)}
,
\qquad
\tilde{\mathcal{D}}_{k,0}^{(s)}=\tilde{\mathcal{B}}_{k,1}^{(s)},
\qquad
\tilde{\mathcal{E}}_{k,0}^{(s)}
=\tilde{\mathcal{B}}_{0,k}^{(s-1)}=\tilde{\mathcal{F}}_{0,k}^{(s-6)}=\tilde{\mathcal{F}}_{k,0}^{(s)}.
\label{trivial relations tilde 1}
\end{gather}
We also have $\mathcal{D}_{0, k}^{(s)} = \tilde{\mathcal{B}}_{k, 1}^{(-s-6k-2)}$, $\tilde{\mathcal{D}}_{0,
k}^{(s)} = \mathcal{B}_{k, 1}^{(-s-6k-2)}$, $k \in \mathbb{Z}_{\geq 0}$, $s \in \mathbb{Z}$.

Note that $\tilde{\mathcal{B}}_{k, \ell}^{(s)}$, $\tilde{\mathcal{D}}_{0, \ell}^{(s)}$,
$\tilde{\mathcal{D}}_{k, 0}^{(s)}$ are minimal af\/f\/inizations.
In general, the modules $\tilde{\mathcal{B}}_{k, \ell}^{(s)}$, $\tilde{\mathcal{C}}_{k, \ell}^{(s)}$,
$\tilde{\mathcal{D}}_{k, \ell}^{(s)}$, $\tilde{\mathcal{E}}_{k, \ell}^{(s)}$, $\tilde{\mathcal{F}}_{k,
\ell}^{(s)}$ are not special.
For example, we have the following proposition.
\begin{Proposition}
The module $\tilde{\mathcal{B}}_{3, 1}^{(0)}=L(1_{0}1_{2}1_{4}2_{11})$ is not special.
\end{Proposition}
\begin{proof}\sloppy
Suppose that $L(1_{0}1_{2}1_{4}2_{11})$ is special.
Then the FM algorithm applies to $L(1_{0}1_{2}1_{4}2_{11})$.
Therefore, by the FM algorithm, the monomials
\begin{gather*}
1_01_21_42_{11},
\quad
1_01_21_6^{-1}2_52_{11},
\quad
1_01_4^{-1}1_6^{-1}2_{3}2_{5}2_{11},
\quad
1_2^{-1}1_4^{-1}1_6^{-1}
2_{1}2_{3}2_{5}2_{11},
\quad
2_7^{-1}2_32_52_{11},
\\
2_{7}^{-1}2_{9}^{-1}1_{4}1_61_82_52_{11},
\quad
2_{7}^{-1}1_41_61_{10}^{-1}2_52_{11},
\quad
1_41_8^{-1}1_{10}^{-1}
2_52_{11},
\quad
1_6^{-1}1_8^{-1}1_{10}^{-1}2_5^{2}2_{11},
\quad
2_{5}
\end{gather*}
are in $\mathscr{M}(L(1_{0}1_{2}1_{4}2_{11}))$.
Hence $\mathscr{M}(L(1_{0}1_{2}1_{4}2_{11}))$ has at least two dominant monomials $1_01_21_42_{11}$ and
$2_{5}$.
This contradicts the assumption that $L(1_{0}1_{2}1_{4}2_{11})$ is special.
\end{proof}
\begin{Theorem}
\label{anti-special}
The modules $\mathcal{\tilde{B}}_{k, \ell}^{(s)}$, $\mathcal{\tilde{C}}_{k, \ell}^{(s)}$,
$\mathcal{\tilde{D}}_{k, \ell}^{(s)}$,
$\mathcal{\tilde{E}}_{k, \ell}^{(s)}$, $\mathcal{\tilde{F}}_{k,\ell}^{(s)}$,
$s\in\mathbb{Z}$, $k,l,\in\mathbb{Z}_{\geq 0}$ are anti-special.
\end{Theorem}
\begin{proof}
This theorem can be proved using the dual arguments of the proof of Theorem~\ref{special}.
\end{proof}
\begin{Lemma}
\label{involution}
Let $\iota: \mathbb{Z}\mathcal{P} \to \mathbb{Z}\mathcal{P}$ be a~homomorphism of rings such that $Y_{1,
aq^{s}} \mapsto Y_{1, aq^{12-s}}^{-1}$, $Y_{2, aq^{s}} \mapsto Y_{2, aq^{12-s}}^{-1}$ for all $a \in
\mathbb{C}^{\times}$, $s\in \mathbb{Z}$.
Then
\begin{gather*}
\chi_q\Big(\tilde{\mathcal{B}}_{k,\ell}^{(s)}\Big)=\iota\Big(\chi_q\Big(\mathcal{B}_{k,\ell}^{(s)}
\Big)\Big),
\qquad
\chi_q\Big(\tilde{\mathcal{C}}_{k,\ell}^{(s)}\Big)=\iota\Big(\chi_q\Big(\mathcal{C}_{k,\ell}^{(s)}\Big)\Big),
\\
\chi_q\Big(\tilde{\mathcal{D}}_{k,\ell}^{(s)}\Big)=\iota\Big(\chi_q\Big(\mathcal{D}_{k,\ell}^{(s)}
\Big)\Big),
\qquad
\chi_q\Big(\tilde{\mathcal{E}}_{k,\ell}^{(s)}\Big)=\iota\Big(\chi_q\Big(\mathcal{E}_{k,\ell}^{(s)}
\Big)\Big),
\qquad
\chi_q\Big(\tilde{\mathcal{F}}_{k,\ell}^{(s)}\Big)=\iota\Big(\chi_q\Big(\mathcal{F}_{k,\ell}^{(s)}\Big)\Big).
\end{gather*}
\end{Lemma}

\begin{proof}
Let $m_+$ be one of $B_{k, \ell}^{(s)}, C_{k, \ell}^{(s)}, D_{k, \ell}^{(s)}, E_{k, \ell}^{(s)}, F_{k,
\ell}^{(s)}$.
Then $\chi_q(\tilde{m_+})$ can be computed by the FM algorithm starting from the lowest weight using $A_{i,
a}$ with $i\in I$, $a\in \mathbb{C}^{\times}$.
The procedure is dual to the computation of $\chi_q(m_+)$ which starts from $m_+$ using $A_{i, a}^{-1}$
with $i\in I$, $a\in \mathbb{C}^{\times}$.
The highest (resp.\
lowest) $l$-weight in $\chi_q(m_+)$ is sent to the lowest (resp.\
highest) $l$-weight in $\chi_q(\tilde{m_+})$ by~$\iota$.
\end{proof}
Note that Lemma~\ref{involution} can also proved using the Cartan involution in~\cite{CP91}.

The modules $\mathcal{\tilde{B}}_{k, \ell}^{(s)}$, $\mathcal{\tilde{C}}_{k, \ell}^{(s)}$,
$\mathcal{\tilde{D}}_{k, \ell}^{(s)}$, $\mathcal{\tilde{E}}_{k, \ell}^{(s)}$, $\mathcal{\tilde{F}}_{k,\ell}^{(s)}$
satisfy the same relations as in Theorem~\ref{extended T-system} but the roles of left and
right modules are exchanged.
More precisely, we have the following theorem.
\begin{Theorem}\label{extended T-system for tilde}
For $s\in\mathbb{Z}$, $k,\ell\in \mathbb{Z}_{\geq 1}$, $t\in\mathbb{Z}_{\geq 2}$, we have the following
relations in $\rep(U_q\hat{\mathfrak{g}})$.
\begin{gather*}
\Big[\tilde{\mathcal{B}}_{k-1,\ell}^{(s+6)}\Big]\Big[\tilde{\mathcal{B}}_{k,\ell-1}^{(s)}\Big]
=\Big[\tilde{\mathcal{B}}_{k,\ell}^{(s)}\Big]\Big[\tilde{\mathcal{B}}_{k-1,\ell-1}^{(s+6)}\Big]
+\Big[\tilde{\mathcal{E}}_{3k-1,\left\lceil\frac{2\ell-2}{3}\right\rceil}^{(s+1)}\Big]
\Big[\tilde{\mathcal{B}}_{\left\lfloor\frac{\ell-1}{3}\right\rfloor,0}^{(s+6k+6)}\Big],
\\
\Big[\tilde{\mathcal{E}}_{0,\ell}^{(s)}\Big]
=\Big[\tilde{\mathcal{B}}_{\left\lfloor\frac{\ell+1}{2}\right\rfloor,0}^{(s+3)}\Big]
\Big[\tilde{\mathcal{B}}_{\left\lfloor\frac{\ell}{2}\right\rfloor,0}^{(s+5)}\Big],
\qquad
\Big[\tilde{\mathcal{E}}_{1, \ell}^{(s)}\Big]
=\Big[\tilde{\mathcal{D}}_{0,\left\lfloor\frac{\ell}{2}\right\rfloor}^{(s-1)}\Big]
\Big[\tilde{\mathcal{B}}_{\left\lfloor\frac{\ell+1}{2}
\right\rfloor,0}^{(s+5)}\Big],
\\
\Big[\tilde{\mathcal{E}}_{t-1,\ell}^{(s+2)}\Big]\Big[\tilde{\mathcal{E}}_{t,\ell-1}^{(s)}\Big]
=\Big[\tilde{\mathcal{E}}_{t,\ell}^{(s)}\Big]\Big[\tilde{\mathcal{E}}_{t-1,\ell-1}^{(s+2)}\Big]
\\
\hphantom{\Big[\tilde{\mathcal{E}}_{t-1,\ell}^{(s+2)}\Big]\Big[\tilde{\mathcal{E}}_{t,\ell-1}^{(s)}\Big]=}{}
+
\begin{cases}
\Big[\tilde{\mathcal{D}}_{r,p-1}^{(s+1)}\Big]\Big[\tilde{\mathcal{B}}_{r+p,0}^{(s+3)}\Big]
\Big[\tilde{\mathcal{B}}_{r,3p-2}^{(s+5)}\Big],
&\text{if}\quad t=3r+2,
\;
\ell=2p-1,
\\[2mm]
\Big[ \tilde{\mathcal{B}}_{r+p+1, 0}^{(s+1)} \Big]
\Big[\tilde{\mathcal{C}}_{r,p}^{(s+3)}\Big]\Big[\tilde{\mathcal{B}}_{r,3p-1}^{(s+5)}\Big],
&\text{if}\quad t=3r+2,
\;
\ell=2p,
\\[2mm]
\Big[ \tilde{\mathcal{B}}_{r+1, 3p-2}^{(s+1)} \Big]
\Big[\tilde{\mathcal{D}}_{r,p-1}^{(s+3)}\Big]\Big[\tilde{\mathcal{B}}_{r+p,0}^{(s+5)}\Big],
&\text{if}\quad t=3r+3,
\;
\ell=2p-1,
\\[2mm]
\Big[ \tilde{\mathcal{B}}_{r+1, 3p-1}^{(s+1)} \Big]
\Big[\tilde{\mathcal{B}}_{r+p+1,0}^{(s+3)}\Big]\Big[\tilde{\mathcal{C}}_{r,p}^{(s+5)}\Big],
&\text{if}\quad t=3r+3,
\;
\ell=2p,
\\[2mm]
\Big[ \tilde{\mathcal{B}}_{r+p+1, 0}^{(s+1)} \Big]
\Big[\tilde{\mathcal{B}}_{r+1,3p-2}^{(s+3)}\Big]\Big[\tilde{\mathcal{D}}_{r,p-1}^{(s+5)}\Big],
&\text{if}\quad t=3r+4,
\;
\ell=2p-1,
\\[2mm]
\Big[ \tilde{\mathcal{C}}_{r+1, p}^{(s+1)} \Big]
\Big[\tilde{\mathcal{B}}_{r+1,3p-1}^{(s+3)}\Big]\Big[\tilde{\mathcal{B}}_{r+p+1,0}^{(s+5)}\Big],
&\text{if}\quad t=3r+4,
\;
\ell=2p,
\end{cases}
\\
\Big[\tilde{\mathcal{C}}_{k-1,\ell}^{(s+6)}\Big]\Big[\tilde{\mathcal{C}}_{k,\ell-1}^{(s)}\Big]
=\Big[\tilde{\mathcal{C}}_{k,\ell}
^{(s)}\Big]\Big[\tilde{\mathcal{C}}_{k-1,\ell-1}^{(s+6)}\Big]+\Big[\tilde{\mathcal{F}}_{3k-2,3\ell-2}^{(s+1)}\Big],
\\
\Big[\tilde{\mathcal{B}}_{\ell,0}^{(s+8)}\Big]\Big[\tilde{\mathcal{D}}_{0,\ell-1}^{(s)}\Big]
=\Big[\tilde{\mathcal{D}}_{0,\ell}
^{(s)}\Big]\Big[\tilde{\mathcal{B}}_{\ell-1,0}^{(s+8)}\Big]+\Big[\tilde{\mathcal{B}}_{0,3\ell-1}^{(s+4)}\Big],
\\
\Big[\tilde{\mathcal{D}}_{k-1,\ell}^{(s+6)}\Big]\Big[\tilde{\mathcal{D}}_{k,\ell-1}^{(s)}\Big]
=\Big[\tilde{\mathcal{D}}_{k,\ell}
^{(s)}\Big]\Big[\tilde{\mathcal{D}}_{k-1,\ell-1}^{(s+6)}\Big]+\Big[\tilde{\mathcal{F}}_{3k-1,3\ell-1}^{(s+1)}\Big],
\\
\Big[\tilde{\mathcal{F}}_{k-1,\ell}^{(s+2)}\Big]\Big[\tilde{\mathcal{F}}_{k,\ell-1}^{(s)}\Big]
=\Big[\tilde{\mathcal{F}}_{k,\ell}
^{(s)}\Big]\Big[\tilde{\mathcal{F}}_{k-1,\ell-1}^{(s+2)}\Big]
\\
\hphantom{\Big[\tilde{\mathcal{F}}_{k-1,\ell}^{(s+2)}\Big]\Big[\tilde{\mathcal{F}}_{k,\ell-1}^{(s)}\Big]=}{}
+
\begin{cases}
\Big[\tilde{\mathcal{B}}_{r,0}^{(s+1)}\Big]\Big[\tilde{\mathcal{D}}_{r,\left\lfloor\frac{\ell}{3}\right\rfloor}^{(s+3)}
\Big]\Big[\tilde{\mathcal{C}}_{r,\left\lfloor\frac{\ell+1}{3}\right\rfloor}^{(s+5)}\Big]
\Big[\tilde{\mathcal{B}}_{\left\lfloor\frac{\ell-1}{3}\right\rfloor,0}^{(s+2k+11)}\Big],
&\text{if}\quad k=3r+1,
\\[2mm]
\Big[ \tilde{\mathcal{C}}_{r+1, \left\lfloor \frac{\ell+1}{3} \right\rfloor}^{(s+1)} \Big]
\Big[\tilde{\mathcal{B}}_{r,0}^{(s+3)}\Big]\Big[\tilde{\mathcal{D}}_{r,\left\lfloor\frac{\ell}{3}\right\rfloor}^{(s+5)}
\Big]\Big[\tilde{\mathcal{B}}_{\left\lfloor\frac{\ell-1}{3}\right\rfloor,0}^{(s+2k+11)}\Big],
&\text{if}\quad k=3r+2,
\\[2mm]
\Big[ \tilde{\mathcal{D}}_{r+1, \left\lfloor \frac{\ell}{3} \right\rfloor}^{(s+1)} \Big]
\Big[\tilde{\mathcal{C}}_{r+1,\left\lfloor\frac{\ell+1}{3}\right\rfloor}^{(s+3)}\Big]
\Big[\tilde{\mathcal{B}}_{r,0}^{(s+5)}\Big]
\Big[\tilde{\mathcal{B}}_{\left\lfloor\frac{\ell-1}{3}\right\rfloor,0}^{(s+2k+11)}\Big],
&\text{if}\quad k=3r+3.
\end{cases}
\end{gather*}
Moreover, the modules corresponding to each summand on the right hand side of the above relations are all
irreducible.
\end{Theorem}
\begin{proof}
The theorem follows from the relations in Theorem~\ref{extended T-system}, Theorem~\ref{irreducible}, and
Lem-\linebreak
ma~\ref{involution}.
\end{proof}

The following proposition is similar to Proposition~\ref{compute}.
\begin{Proposition}
Given $\chi_q(1_s)$, $\chi_q(2_s)$, one can obtain the $q$-characters of $\tilde{\mathcal{B}}_{k,
\ell}^{(s)}$, $\tilde{\mathcal{C}}_{k, \ell}^{(s)}$, $\tilde{\mathcal{D}}_{k, \ell}^{(s)}$,
$\tilde{\mathcal{E}}_{k, \ell}^{(s)}$, $\tilde{\mathcal{F}}_{k, \ell}^{(s)}$,
$s\in\mathbb{Z}$, $k,\ell\in \mathbb{Z}_{\geq 0}$, recursively, by using~\eqref{trivial relations tilde 1},
and computing the
$q$-character of the top module through the $q$-characters of other modules in relations in
Theorem~{\rm \ref{extended T-system for tilde}}.
\end{Proposition}

\section{Dimensions}
\label{dimensions section}

In this section, we give dimension formulas for the modules $\mathcal{B}_{k, \ell}^{(s)}$, $\mathcal{C}_{k,
\ell}^{(s)}$, $\mathcal{D}_{k, \ell}^{(s)}$, $\mathcal{E}_{k, \ell}^{(s)}$, $\mathcal{F}_{k, \ell}^{(s)}$,
$\tilde{\mathcal{B}}_{k, \ell}^{(s)}$, $\tilde{\mathcal{C}}_{k, \ell}^{(s)}$, $\tilde{\mathcal{D}}_{k,
\ell}^{(s)}$, $\tilde{\mathcal{E}}_{k, \ell}^{(s)}$, $\tilde{\mathcal{F}}_{k, \ell}^{(s)}$.

Note that dimensions do not depend on the upper index $s$.
Note also that $\dim M = \dim \tilde{M}$ for each $M=\mathcal{B}_{k, \ell}^{(s)}$, $\mathcal{C}_{k,
\ell}^{(s)}$, $\mathcal{D}_{k, \ell}^{(s)}$, $\mathcal{E}_{k, \ell}^{(s)}$, $\mathcal{F}_{k, \ell}^{(s)}$.
\begin{Theorem}
\label{dimensions}
Let $s\in\mathbb{Z}$, $k,\ell\in \mathbb{Z}_{\geq 0}$.
Then
\begin{gather*}
\dim\mathcal{B}_{k,3\ell}^{(s)}=(\ell+2)(\ell+1)(1+k)(k+3+\ell)(k+2+\ell)
\\
\hphantom{\dim\mathcal{B}_{k,3\ell}^{(s)}=}{}
\times\big(54\ell^3k^3+243\ell^2k^3+363\ell k^3+180k^3+2784\ell^2k^2+1080k^2+162\ell^4k^2
\\
\hphantom{\dim\mathcal{B}_{k,3\ell}^{(s)}=}{}
+2880\ell k^2+1134\ell^3k^2+162\ell^5k+1539\ell^4k+5490\ell^3k+9132\ell^2k+7057\ell k
\\
\hphantom{\dim\mathcal{B}_{k,3\ell}^{(s)}=}{}
+2040k+54\ell^6+648\ell^5+3069\ell^4+7272\ell^3+8977\ell^2+5380\ell+1200\big)/14400,
\\
\dim\mathcal{B}_{k,3\ell+1}^{(s)}=(\ell+3)(\ell+2)(\ell+1)(1+k)(k+2+\ell)(k+4+\ell)(k+3+\ell)
\\
\hphantom{\dim\mathcal{B}_{k,3\ell+1}^{(s)}=}{}
\times\big(171\ell k^2+120k^2+54\ell^2k^2+600k+621\ell^2k+108\ell^3k
\\
\hphantom{\dim\mathcal{B}_{k,3\ell+1}^{(s)}=}{}
+1116\ell k+54\ell^4+450\ell^3+1341\ell^2+1665\ell+700\big)/14400,
\\
\dim\mathcal{B}_{k,3\ell+2}^{(s)}=(\ell+3)(\ell+2)(\ell+1)(1+k)(k+4+\ell)(k+3+\ell)(2+k+\ell)
\\
\hphantom{\dim\mathcal{B}_{k,3\ell+2}^{(s)}=}{}
\times\big(300k^2+261\ell k^2+54\ell^2k^2+891\ell^2k+2376\ell k+2040k
\\
\hphantom{\dim\mathcal{B}_{k,3\ell+2}^{(s)}=}{}
+108\ell^3k+54\ell^4+630\ell^3+2691\ell^2+4995\ell+3400\big)/14400,
\\
\dim\mathcal{C}_{k,\ell}^{(s)}=(\ell+2)(\ell+1)(k+2)(k+1)(k+3+\ell)(k+2+\ell)
\\
\hphantom{\dim\mathcal{C}_{k,\ell}^{(s)}=}{}
\times \big(3k^2+3\ell k^2+12k+15\ell k+3\ell^2k+3\ell^2+12\ell+10\big)/240,
\\
\dim\mathcal{D}_{k,\ell}^{(s)}=(\ell+2)(\ell+1)(k+2)(k+1)(k+3+\ell)(k+4+\ell)
\\
\hphantom{\dim\mathcal{D}_{k,\ell}^{(s)}=}{}
\times \big(3\ell k^2+6k^2+3\ell^2k+30k+21\ell k+6\ell^2+30\ell+35\big)/240,
\\
\dim\mathcal{E}_{3k,2\ell}^{(s)}=(\ell+2)(\ell+1)(k+1)(k+\ell+1)(k+\ell+2)^2(k+\ell+3)^2
\\
\hphantom{\dim\mathcal{E}_{3k,2\ell}^{(s)}=}{}
\times\big(27k^4\ell^2+81k^4\ell+54k^4+81k^3\ell^3+468k^3\ell^2+825k^3\ell
\\
\hphantom{\dim\mathcal{E}_{3k,2\ell}^{(s)}=}{}
+432k^3+81k^2\ell^4+711k^2\ell^3+2184k^2\ell^2+2754k^2\ell+1179k^2+27k\ell^5
\\
\hphantom{\dim\mathcal{E}_{3k,2\ell}^{(s)}=}{}
+342k\ell^4+1593k\ell^3+3438k\ell^2+3435k\ell+1260k+18\ell^5
\\
\hphantom{\dim\mathcal{E}_{3k,2\ell}^{(s)}=}{}
+180\ell^4+696\ell^3+1296\ell^2+1160\ell+400\big)/28800,
\\
\dim\mathcal{E}_{3k,2\ell+1}^{(s)}=(\ell+3)(\ell+2)(\ell+1)(k+1)(k+\ell+4)(k+\ell+2)^2(k+\ell+3)^2
\\
\hphantom{\dim\mathcal{E}_{3k,2\ell+1}^{(s)}=}{}
\times\big(27k^4\ell+54k^4+81k^3\ell^2+414k^3\ell+510k^3+81k^2\ell^3
\\
\hphantom{\dim\mathcal{E}_{3k,2\ell+1}^{(s)}=}{}
+684k^2\ell^2+1842k^2\ell+1611k^2+27k\ell^4+342k\ell^3
\\
\hphantom{\dim\mathcal{E}_{3k,2\ell+1}^{(s)}=}{}
+1512k\ell^2+2808k\ell+1875k+18\ell^4+180\ell^3+642\ell^2+960\ell+500\big)/28800,
\\
\dim\mathcal{E}_{3k+1,2\ell}^{(s)}=(\ell+2)(\ell+1)(k+1)(k+\ell+4)(k+\ell+2)^2(k+\ell+3)^2
\\
\hphantom{\dim\mathcal{E}_{3k+1,2\ell}^{(s)}=}{}
\times\big(27k^4\ell^2+81k^4\ell+54k^4+81k^3\ell^3+477k^3\ell^2+852k^3\ell
\\
\hphantom{\dim\mathcal{E}_{3k+1,2\ell}^{(s)}=}{}
+450k^3+81k^2\ell^4+747k^2\ell^3+2373k^2\ell^2+3069k^2\ell+1341k^2
\\
\hphantom{\dim\mathcal{E}_{3k+1,2\ell}^{(s)}=}{}
+27k\ell^5+387k\ell^4+1935k\ell^3+4353k\ell^2+4461k\ell+1665k
\\
\hphantom{\dim\mathcal{E}_{3k+1,2\ell}^{(s)}=}{}
+36\ell^5+360\ell^4+1374\ell^3+2490\ell^2+2140\ell+700\big)/28800,
\\
\dim\mathcal{E}_{3k+1,2\ell+1}^{(s)}=(\ell+3)(\ell+2)(\ell+1)(k+1)(k+\ell+2)(k+\ell+3)^2(k+\ell+4)^2
\\
\hphantom{\dim\mathcal{E}_{3k+1,2\ell+1}^{(s)}=}{}
\times\big(27k^4\ell+54k^4+81k^3\ell^2+450k^3\ell+582k^3+81k^2\ell^3+774k^2\ell^2
\\
\hphantom{\dim\mathcal{E}_{3k+1,2\ell+1}^{(s)}=}{}
+2310k^2\ell+2193k^2+27k\ell^4+414k\ell^3+2124k\ell^2+4488k\ell
\\
\hphantom{\dim\mathcal{E}_{3k+1,2\ell+1}^{(s)}=}{}
+3375k+36\ell^4+396\ell^3+1590\ell^2+2760\ell+1750\big)/28800,
\\
\dim\mathcal{E}_{3k+2,2\ell}^{(s)}=(\ell+2)(\ell+1)(k+2)(k+1)(k+\ell+4)(k+\ell+2)(k+\ell+3)^2
\\
\hphantom{\dim\mathcal{E}_{3k+2,2\ell}^{(s)}=}{}
\times\big(27k^4\ell^2+81k^4\ell+54k^4+108k^3\ell^3+648k^3\ell^2+1176k^3\ell
\\
\hphantom{\dim\mathcal{E}_{3k+2,2\ell}^{(s)}=}{}
+630k^3+162k^2\ell^4+1458k^2\ell^3+4629k^2\ell^2+6057k^2\ell
\\
\hphantom{\dim\mathcal{E}_{3k+2,2\ell}^{(s)}=}{}
+
2691k^2+108k\ell^5+1296k\ell^4+5946k\ell^3+12942k\ell^2+13230k\ell+4995k
\\
\hphantom{\dim\mathcal{E}_{3k+2,2\ell}^{(s)}=}{}
+27\ell^6+405\ell^5+2439\ell^4+7515\ell^3+12429\ell^2+10395\ell+3400\big)/28800,
\\
\dim\mathcal{E}_{3k+2,2\ell+1}^{(s)}=(\ell+3)(\ell+2)(\ell+1)(k+2)(k+1)(k+\ell+5)(k+\ell+2)
\\
\hphantom{\dim\mathcal{E}_{3k+2,2\ell+1}^{(s)}=}{}
\times\big(k+\ell+3)^2(k+\ell+4)^2(9k^2\ell+18k^2
\\
\hphantom{\dim\mathcal{E}_{3k+2,2\ell+1}^{(s)}=}{}
+18k\ell^2+99k\ell+128k+9\ell^3+81\ell^2+237\ell+225\big)/9600,
\\
\dim\mathcal{F}_{3k,3\ell}^{(s)}=(\ell+2)^2(\ell+1)^2(k+2)^2(k+1)^2(k+\ell+3)^2
\\
\hphantom{\dim\mathcal{F}_{3k,3\ell}^{(s)}=}{}
\times\big(27k^4\ell^2+81k^4\ell+54k^4+54k^3\ell^3+405k^3\ell^2+801k^3\ell+432k^3+27k^2\ell^4
\\
\hphantom{\dim\mathcal{F}_{3k,3\ell}^{(s)}=}{}
+405k^2\ell^3+1746k^2\ell^2+2646k^2\ell+1179k^2+81k\ell^4+801k\ell^3
\\
\hphantom{\dim\mathcal{F}_{3k,3\ell}^{(s)}=}{}
+2646k\ell^2+3342k\ell+1260k+54\ell^4+432\ell^3+1179\ell^2+1260\ell+400\big)/57600,
\\
\dim\mathcal{F}_{3k+1,3\ell}^{(s)}=(\ell+2)^2(\ell+1)^2(k+3)(k+1)(k+2)^2(k+\ell+4)(k+\ell+3)
\\
\hphantom{\dim\mathcal{F}_{3k+1,3\ell}^{(s)}=}{}
\times\big(27k^4\ell^2+81k^4\ell+54k^4+54k^3\ell^3+414k^3\ell^2+828k^3\ell+450k^3+
27k^2\ell^4
\\
\hphantom{\dim\mathcal{F}_{3k+1,3\ell}^{(s)}=}{}
+414k^2\ell^3+1854k^2\ell^2+2907k^2\ell+1341k^2+81k\ell^4+864k\ell^3
+
3063k\ell^2\\
\hphantom{\dim\mathcal{F}_{3k+1,3\ell}^{(s)}=}{}
+4116k\ell+1665k+54\ell^4+498\ell^3+1563\ell^2+1905\ell+700\big)/57600,
\\
\dim\mathcal{F}_{3k+2,3\ell}^{(s)}=(\ell+2)^2(\ell+1)^2(k+3)(k+1)(k+2)^2(k+\ell+4)(k+\ell+3)
\\
\hphantom{\dim\mathcal{F}_{3k+2,3\ell}^{(s)}=}{}
\times\big(27k^4\ell^2+81k^4\ell+54k^4+54k^3\ell^3+504k^3\ell^2+1098k^3\ell+630k^3+27k^2\ell^4
\\
\hphantom{\dim\mathcal{F}_{3k+2,3\ell}^{(s)}=}{}
+
558k^2\ell^3+3042k^2\ell^2+5355k^2\ell+2691k^2+135k\ell^4+1764k\ell^3+7395k\ell^2
\\
\hphantom{\dim\mathcal{F}_{3k+2,3\ell}^{(s)}=}{}
+
11190k\ell+4995k+162\ell^4+1734\ell^3+6249\ell^2+8475\ell+3400\big)/57600,
\\
\dim\mathcal{F}_{3k,3\ell+1}^{(s)}=(\ell+3)(\ell+1)(\ell+2)^2(k+2)^2(k+1)^2(k+\ell+4)(k+\ell+3)
\\
\hphantom{\dim\mathcal{F}_{3k,3\ell+1}^{(s)}=}{}
\times\big(27k^4\ell^2+81k^4\ell+54k^4+54k^3\ell^3+414k^3\ell^2+864k^3\ell+498k^3
+
27k^2\ell^4
\\
\hphantom{\dim\mathcal{F}_{3k,3\ell+1}^{(s)}=}{}
+414k^2\ell^3+1854k^2\ell^2+3063k^2\ell+1563k^2+81k\ell^4+828k\ell^3
+
2907k\ell^2\\
\hphantom{\dim\mathcal{F}_{3k,3\ell+1}^{(s)}=}{}
+4116k\ell+1905k+54\ell^4+450\ell^3+1341\ell^2+1665\ell+700\big)/57600,
\\
\dim\mathcal{F}_{3k+1,3\ell+1}^{(s)}=(\ell+3)(\ell+1)(\ell+2)^2(k+3)(k+1)(k+2)^2(k+\ell+3)(k+\ell+4)
\\
\hphantom{\dim\mathcal{F}_{3k+1,3\ell+1}^{(s)}=}{}
\times\big(27k^4\ell^2\!+81k^4\ell+54k^4+54k^3\ell^3\!+450k^3\ell^2+972k^3\ell+570k^3\!+27k^2\ell^4
\\
\hphantom{\dim\mathcal{F}_{3k+1,3\ell+1}^{(s)}=}{}
+
450k^2\ell^3+2214k^2\ell^2+3891k^2\ell+2061k^2+81k\ell^4+972k\ell^3+3891k\ell^2
\\
\hphantom{\dim\mathcal{F}_{3k+1,3\ell+1}^{(s)}=}{}
+
6060k\ell+2985k+54\ell^4+570\ell^3+2061\ell^2+2985\ell+1400\big)/57600,
\\
\dim\mathcal{F}_{3k+2,3\ell+1}^{(s)}=(\ell+3)(\ell+1)(\ell+2)^2(k+3)(k+1)(k+2)^2(k+\ell+3)(k+\ell+5)
\\
\hphantom{\dim\mathcal{F}_{3k+2,3\ell+1}^{(s)}=}{}
\times\big(k+\ell+4)^2(9k^2\ell^2+27k^2\ell+18k^2+45k\ell^2+135k\ell
\\
\hphantom{\dim\mathcal{F}_{3k+2,3\ell+1}^{(s)}=}{}
+88k+54\ell^2+164\ell+105\big)/19200,
\\
\dim\mathcal{F}_{3k,3\ell+2}^{(s)}=(\ell+3)(\ell+1)(\ell+2)^2(k+2)^2(k+1)^2(k+\ell+4)(k+\ell+3)\big(27k^4\ell^2
\\
\hphantom{\dim\mathcal{F}_{3k,3\ell+2}^{(s)}=}{}
 +135k^4\ell+162k^4+54k^3\ell^3+558k^3\ell^2+1764k^3\ell+1734k^3+27k^2\ell^4
\\
\hphantom{\dim\mathcal{F}_{3k,3\ell+2}^{(s)}=}{}
+504k^2\ell^3+3042k^2\ell^2+7395k^2\ell+6249k^2+81k\ell^4+1098k\ell^3+5355k\ell^2
\\
\hphantom{\dim\mathcal{F}_{3k,3\ell+2}^{(s)}=}{}
+11190k\ell+8475k+54\ell^4+630\ell^3+2691\ell^2+4995\ell+3400\big)/57600,
\\
\dim\mathcal{F}_{3k+1,3\ell+2}^{(s)}=(\ell+3)(\ell+1)(\ell+2)^2(k+3)(k+1)(k+2)^2(k+\ell+3)(k+\ell+5)
\\
\hphantom{\dim\mathcal{F}_{3k+1,3\ell+2}^{(s)}=}{}
\times\big(k+\ell+4)^2(9k^2\ell^2+45k^2\ell+54k^2+27k\ell^2+135k\ell
\\
\hphantom{\dim\mathcal{F}_{3k+1,3\ell+2}^{(s)}=}{}
+164k+18\ell^2+88\ell+105\big)/19200,
\\
\dim\mathcal{F}_{3k+2,3\ell+2}^{(s)}=(\ell+3)(\ell+1)(\ell+2)^2(k+3)(k+1)(k+2)^2(k+\ell+4)(k+\ell+5)
\\
\hphantom{\dim\mathcal{F}_{3k+2,3\ell+2}^{(s)}=}{}
\times\big(27k^4\ell^2+135k^4\ell+162k^4+54k^3\ell^3+630k^3\ell^2+2124k^3\ell+2166k^3
\\
\hphantom{\dim\mathcal{F}_{3k+2,3\ell+2}^{(s)}=}{}
+27k^2\ell^4+630k^2\ell^3+4374k^2\ell^2+11661k^2\ell+10473k^2+135k\ell^4
\\
\hphantom{\dim\mathcal{F}_{3k+2,3\ell+2}^{(s)}=}{}
+2124k\ell^3+11661k\ell^2+26748k\ell+21759k+162\ell^4+2166\ell^3
\\
\hphantom{\dim\mathcal{F}_{3k+2,3\ell+2}^{(s)}=}{}
+10473\ell^2+21759\ell+16400\big)/57600.
\end{gather*}
\end{Theorem}
\begin{proof}
We check the initial conditions, namely dimensions of $\mathcal{B}^{(s)}_{0,1}$, $\mathcal{B}^{(s)}_{1,0}$.
We check the dimensions are compatible with relations~\eqref{trivial relations 1},~\eqref{T-sys
1},~\eqref{T-sys 2}.
We directly check that the formulas satisfy the relations in Theorems~\ref{extended T-system}.
For the checks we employed the computer algebra system Maple.

The theorem follows since the solution of the extended $T$-system is unique, see
Proposi\-tion~\ref{compute}.
\end{proof}

\subsection*{Acknowledgements}

We would like to thank D.~Hernandez, B.~Leclerc, T.~Nakanishi, C.A.S.~Young for helpful discussions.
JL would like to thank IUPUI Department of Mathematical Sciences for hospitality during his visit when this
work was carried out.
JL is partially supported by a~CSC scholarship, the Natural Science Foundation of Gansu Province
(No.~1107RJZA218), and the Fundamental Research Funds for the Central Universities (No.~lzujbky-2012-12)
from China.
The research of EM is supported by the NSF, grant number DMS-0900984.

\newpage

\pdfbookmark[1]{References}{ref}
\LastPageEnding

\end{document}